\documentclass{amsart}

\usepackage{pdfsync}
\newcommand{\s}{\vspace{0.3cm}}
\newtheorem{theo}{Theorem}[section]

\newtheorem{rema}[theo]{Remark}

\numberwithin{equation}{section}

\begin{document}

\title{A computational note about Fricke-Macbeath's curve}
\author{Ruben A. Hidalgo}


\address{Departamento de Matem\'atica, Universidad T\'ecnica Federico Santa Mar\'{\i}a. Casilla 110-V, Valpara\'{\i}so, Chile}
\email{ruben.hidalgo@usm.cl}

\begin{abstract}
The well known Hurwitz upper bound states that a closed Riemann surface $S$ of genus $g \geq 2$ has at most $84(g-1)$ conformal automorphisms. If $S$ has exactly $84(g-1)$ conformal automorphisms, then it is called a Hurwitz curve. 
The first two genera for which there are Hurwitz's curves are $g \in \{3,7\}$. In both situations there is exactly one such curve up to conformal equivalence, in particular, in both cases the field of moduli is ${\mathbb Q}$. As these two curves are quasiplatonic curves, they are definable over ${\mathbb Q}$. The Hurwitz's curve of genus $g=3$  is given by Klein's quartic $x^3y+y^3z+z^3x=0$. The Hurwitz's curve of genus  $g=7$ is known as Fricke-Macbeath's curve and equations over ${\mathbb Q}(\rho)$, where $\rho=e^{2 \pi i/7}$, are known due to Macbeath. Unfortunately, explicit equations over ${\mathbb Q}$ are not easy to find for this curve. In this paper we first explain how to construct an explicit model $Z_{2}$ of Fricke-Macbeath's curve over ${\mathbb Q}(\sqrt{-7})$ and an explicit isomorphism $L_{1}:X \to Z_{2}$, defined over ${\mathbb Q}(\rho)$. Next, using that explicit model we construct another explicit isomorphism $L_{2}:Z_{2} \to W$, defined over ${\mathbb Q}(\sqrt{-7})$, where $W$ is some algebraic curve defined over ${\mathbb Q}$.  Unfortunately, the equations for $W$ are quite long to write down, but everything is explained in order to perform the computations in a computer.

\end{abstract}

\maketitle

\section{Introduction}

Let $S$ be a closed  Riemann surface of genus $g \geq 2$. It is well known that $|{\rm Aut}(S)| \leq 84(g-1)$ (Hurwitz's upper bound). If happens that  $|{\rm Aut}(S)|=84(g-1)$, then one says that $S$ is a Hurwitz's curve. In this case, $S/{\rm Aut}(S)$ is an orbifold with signature $(0;2,3,7)$, that is, $S={\mathbb H}^{2}/\Gamma$, where $\Gamma$ is a torsion free normal subgroup of finite index in the triangular Fuchsian group $\Delta=\langle x,y: x^2=y^3=(xy)^7=1\rangle$ acting on the hyperbolic plane ${\mathbb H}^{2}$ as isometries.
As already noticed by Wiman \cite{Wiman}, in genera $g=2,4,5,6$ there are no Hurwitz's curves and that for $g=3$ there is exactly one (up to holomorphic equivalence); this being Klein's quartic $x^3y+y^3z+z^3x=0$; whose automorphisms group is the  simple group ${\rm PSL}(2,7)$ (of order $168$).

 In genus $g=7$ there is only one (up to conformal equivalence) Hurwitz's curve, called Fricke-Macbeath's curve \cite{Macbeath}. It follows from this uniqueness that the field of moduli of Fricke-Macbeath's curve is ${\mathbb Q}$, the field of rational numbers. As quasiplatonic curves can be defined over their fields of modui \cite{Wolfart} and Hurwitz's curve are quasiplatonic curves, it follows that Fricke-Macbeath's curve can be defined over ${\mathbb Q}$. It seems that in the literature no 
equations over ${\mathbb Q}$ are written for this curve. The automorphisms group of Fricke-Macbeath's curve is the simple group ${\rm PSL}(2,8)$, consisting of $504$ symmetries.
In \cite{Macbeath} Macbeath computed the following explicit equations over ${\mathbb Q}(\rho)$, where $\rho=e^{2 \pi i/7}$,  for Fricke-Macbeath's curve involving three particular elliptic curves as follows: 
\begin{equation}\label{eq0}
X=\left\{\begin{array}{lcl}
y_{1}^{2}=(x-1)(x-\rho^3)(x-\rho^5)(x-\rho^6)\\
\\
y_{2}^{2}=(x-\rho^{2})(x-\rho^4)(x-\rho^5)(x-\rho^6)\\
\\
y_{4}^{2}=(x-\rho)(x-\rho^3)(x-\rho^4)(x-\rho^5)
\end{array}
\right\} \subset {\mathbb C}^{4}.
\end{equation}

In the talk \cite{Wolfart-pres} there is a misprint for the first of the elliptic curves in equations. In the above model is easy to see a group $G \cong {\mathbb Z}_{2}^{3}$ of holomorphic automorphisms of $X$ generated by   
$A_{1}(x,y_{1},y_{2},y_{4})=(x,-y_{1},y_{2},y_{4})$, 
$A_{2}(x,y_{1},y_{2},y_{4})=(x,y_{1},-y_{2},y_{4})$,
and $A_{3}(x,y_{1},y_{2},y_{4})=(x,y_{1},y_{2},-y_{4})$. 

In Section \ref{appendix} we provide a rough explanation about the elliptic curves in the above equations (different from the approach in \cite{Macbeath}) in geometric terms of the highest regular branched Abelian cover of the orbifold $X/G$ of signature $(0;2,2,2,2,2,2,2)$. 

An automorphism of order $7$ of Fricke-Macbeath's curve is given in such model by
$$B(x,y_{1},y_{2},y_{4})=\left(\rho x,\rho^{2} y_{2},\rho^{2} y_{4},\rho^{2} \frac{y_{1}y_{2}}{(x-\rho^{5}) (x-\rho^{6})}\right).$$

The automorphism $B$ normalizes $G$ and it induces, on the orbifold $X/G=\widehat{\mathbb C}$, the rotation $T(x)=\rho x$; moreover, $X/\langle G,B\rangle$
has signature $(0;2,7,7)$, that is, the group $\langle G,B\rangle$ defines a regular dessin d'enfants $(X,\beta)$, where $\beta(x,y_{1},y_{2},y_{4})=x^{7}$ (called an Edmods map) \cite{Wolfart-ankara,Wolfart-pres}. Such a dessin d'enfants is seen to be defined over ${\mathbb Q}(\rho)$, but it is known to be definable over ${\mathbb Q}(\sqrt{-7})$ \cite{Wolfart-ankara}. As a direct consequence of our computations, we obtain an isomorphic dessin defined over ${\mathbb Q}(\sqrt{-7})$ (see Remark \ref{dibujo} in Section \ref{Sec:Z2}).

Notice that $X$ also admits the following anticonformal involution
$$J(x,y_{1},y_{2},y_{4})=\left(\dfrac{1}{\overline{x}}, \dfrac{\overline{y_{1}}}{\overline{x}^{2}},
\dfrac{\rho^{5} \overline{y_{2}}}{\overline{x}^{2}},\dfrac{\rho^{3} \overline{y_{4}}}{\overline{x}^{2}}\right).$$

We may see that $J B J=B$ and $J A_{j} J=A_{j}$, for $j=1,2,4$.
We also have the regular dessin d'enfants $(X,\delta)$, where $\delta(x,y_{1},y_{2},y_{4})=1/x^{7}$. As $\delta=C \circ \beta \circ J$, where $C(x)=\overline{x}$, we have that these two dessins are chirals (these are the two dessins in genus $7$ whose graph is the complete graph $K_{8}$ appearing in \cite{Wolfart-ankara}).

\s

In this paper  we first describe a theoretical/computational method which permits to obtain an explicit a birational isomorphism, defined over ${\mathbb Q}(\rho)$) $$L:X \to W$$ so that  $W$ satisfies that $W^{\sigma}=W$ for every $\sigma \in {\rm Gal}({\mathbb Q}(\rho)/{\mathbb Q})$; that is, $W$ is defined over ${\mathbb Q}$. This method is based on the constructive proof of Weil's Galois descent theorem \cite{HR}. In order to get explicitly $L$, one needs to find an explicit set of generators of the invariant polynomials under a suitable linear cyclic group of order $6$ acting by permutations on a $24$-dimensional space. Unfortunately, it is hard to get (inclusive with MAGMA \cite{MAGMA}) such a set of generators. Instead to proceed in a direct way, we divide this search into two parts. In the first one we explain how to provide (with the help of a computer) an explicit algebraic curve $Z_{2}$ defined over ${\mathbb Q}(\sqrt{-7})$ together an explicit isomorphism $L_{1}^{*}:X \to Z_{2}$ (also we provide explicitly its inverse $(L_{1}^{*})^{-1}:Z_{2} \to X$). This is the first known model of Fricke-Macbeath's curve over ${\mathbb Q}(\sqrt{-7})$ to our knowledge. Also, the isomorphism $L_{1}^{*}$ provides an isomorphism between the dessin $(X,\beta)$ with the dessin $(Z_{2},\beta^{*})$, where $\beta^{*}$ is defined over ${\mathbb Q}$, that is, the dessin is defined over ${\mathbb Q}(\sqrt{-7})$.
Next, using the explicit model $Z_{2}$, we provide an explicit isomorphism $L_{2}:Z_{2} \to W$, where $W$ is defined over ${\mathbb Q}$ (we are able to see that $W$ is defined over ${\mathbb Q}$ without the knowledge of equations for it). As both, $Z_{2}$ and $L_{2}$ are explicitly given, equations for $W$ over ${\mathbb Q}$ are possible to compute (with the help of a computer), but unfortunately they are long to write them in this paper. Now, $L=L_{2} \circ L_{1}^{*}:X \to W$ is an explicit isomorphism as desired.  

Recently, Bradley Brock told me that he was able to compute a plane equation for Fricke-Macbeath's equation over ${\mathbb Q}$ as
$1 + 7xy + 21x^2 y^2 + 35x^3 y^3 + 28x^4 y^4 + 2x^7 + 2y^7 = 0$. To obtain the model he used the canonical
model in Macbeath's original paper \cite{Macbeath} and then got a suitable change of variables. Moreover, he asserted that the jacobian is isogenous to $E^7$ where $j(E)=1792$
and $E $does not have Complex Multiplication.

\section{In search of a birational isomorphism}
In this section we explain the general ideas in the search of an isomorphism $L:X \to W$, where $W$ is defined over ${\mathbb Q}$. The same ideas will be used in the next sections to obtain it explicitly by working it in two steps.

The Galois extension ${\mathbb Q}(\rho)/{\mathbb Q}$ has as Galois group $\Gamma={\rm Gal}({\mathbb Q}(\rho)/{\mathbb Q})=\langle \sigma \rangle \cong {\mathbb Z}_{6}$, where
$\sigma(\rho)=\rho^{3}$. In particular, $\sigma^{2}(\rho)=\rho^2$, $\sigma^{3}(\rho)=\rho^6$, $\sigma^{4}(\rho)=\rho^4$ and $\sigma^{5}(\rho)=\rho^5$.

\subsection{A weil's datum for $X$}
All the curves $X$, $X^{\sigma}$, $X^{\sigma^{2}}$, $X^{\sigma^{3}}$, $X^{\sigma^{4}}$ and $X^{\sigma^{5}}$ are birationally equivalent. In fact, we may describe explicit birational isomorphisms as follows.

$$f_{\sigma}:X \to X^{\sigma}$$
$$(x,y_{1},y_{2},y_{4}) \mapsto \left( \frac{1}{x}, \frac{y_{1}}{x^{2}}, 
\frac{\rho y_{2}y_{4}}{x^{2}(x-\rho^{4})(x-\rho^{5})}, \frac{\rho^{2} y_{2}}{x^{2}} \right)$$

$$f_{\sigma^{2}}:X \to X^{\sigma^{2}}$$
$$(x,y_{1},y_{2},y_{4}) \mapsto \left( x,y_{1},y_{4},\frac{y_{2}y_{4}}{(x-\rho^{4})(x-\rho^{5})}\right)$$

$$f_{\sigma^{3}}:X \to X^{\sigma^{3}}$$
$$(x,y_{1},y_{2},y_{4}) \mapsto \left( \frac{1}{x}, \frac{y_{1}}{x^{2}}, \frac{\rho^{2} y_{2}}{x^{2}}, 
\frac{\rho^{4} y_{4}}{x^{2}} \right)$$

$$f_{\sigma^{4}}:X \to X^{\sigma^{4}}$$
$$(x,y_{1},y_{2},y_{4}) \mapsto \left( x,y_{1},\frac{y_{2}y_{4}}{(x-\rho^{4})(x-\rho^{5})},y_{2}\right)$$

$$f_{\sigma^{5}}:X \to X^{\sigma^{5}}$$
$$(x,y_{1},y_{2},y_{4}) \mapsto \left( \frac{1}{x}, \frac{y_{1}}{x^{2}}, \frac{\rho^{4} y_{4}}{x^{2}}, 
\frac{\rho y_{2}y_{4}}{x^{2}(x-\rho^{4})(x-\rho^{5})} \right)$$

It can be easily checked that the following Weil's co-cycle condition holds

$$f_{\sigma_{1}\sigma_{2}}=f_{\sigma_{2}}^{\sigma_{1}} \circ f_{\sigma_{1}}, \mbox{ for every }\sigma_{1},\sigma_{2} \in \Gamma.$$

In fact, the last co-cycle condition was used to define $f_{\sigma^{j}}$ starting from $f_{\sigma}$. In particular, $f_{e}=I$, where $e=\sigma^{0}$ is the identity of ${\rm Gal}({\mathbb Q}(\rho)/{\mathbb Q})$ and $I$ is the identity automorphism of $X$. In other words, the collection of birational isomorphisms $\{f_{\sigma^{j}}\}_{j=0}^{5}$ is a Weil's datum for $X$ (this is another way to see that $X$ is definable over ${\mathbb Q}$ as a consequence of Weil's Galois descent theorem \cite{Weil}).

\subsection{Another model for $X$}
Let us consider the rational map

$$\Phi:X \subset {\mathbb C}^{4} \to {\mathbb C}^{24}$$
$$(x,y_{1},y_{2},y_{4}) \mapsto (\vec{x},\vec{z},\vec{w},\vec{u},\vec{v},\vec{r})$$
where
$$\vec{x}=(x_{1},x_{2},x_{3},x_{4})=(x,y_{1},y_{2},y_{4})$$

$$\vec{z}=(z_{1},z_{2},z_{3},z_{4})=f_{\sigma}(\vec{x})=\left( \frac{1}{x}, \frac{y_{1}}{x^{2}}, 
\frac{\rho y_{2}y_{4}}{x^{2}(x-\rho^{4})(x-\rho^{5})}, \frac{\rho^{2} y_{2}}{x^{2}} \right)$$

$$\vec{w}=(w_{1},w_{2},w_{3},w_{4})=f_{\sigma^{2}}(\vec{x})= \left( x,y_{1},y_{4},\frac{y_{2}y_{4}}{(x-\rho^{4})(x-\rho^{5})}\right)$$

$$\vec{u}=(u_{1},u_{2},u_{3},u_{4})=f_{\sigma^{3}}(\vec{x})=\left( \frac{1}{x}, \frac{y_{1}}{x^{2}}, \frac{\rho^{2} y_{2}}{x^{2}}, 
\frac{\rho^{4} y_{4}}{x^{2}} \right)$$

$$\vec{v}=(v_{1},v_{2},v_{3},v_{4})=f_{\sigma^{4}}(\vec{x})=\left( x,y_{1},\frac{y_{2}y_{4}}{(x-\rho^{4})(x-\rho^{5})},y_{2}\right)$$

$$\vec{r}=(r_{1},r_{2},r_{3},r_{4})=f_{\sigma^{5}}(\vec{x})=\left( \frac{1}{x}, \frac{y_{1}}{x^{2}}, \frac{\rho^{4} y_{4}}{x^{2}}, 
\frac{\rho y_{2}y_{4}}{x^{2}(x-\rho^{4})(x-\rho^{5})} \right)$$

It turns out that $\Phi:X \to \Phi(X)$ is a birational isomorphism (the inverse is just given by the projection on the $\vec{x}$-coordinate).

Equations defining the algebraic curve $\Phi(X)$ are the following ones

\begin{equation}\label{eq2}
\Phi(X)=\left\{\begin{array}{lcl}
x_{2}^{2}=(x_{1}-1)(x_{1}-\rho^3)(x_{1}-\rho^5)(x_{1}-\rho^6)\\
\\
x_{3}^{2}=(x_{1}-\rho^{2})(x_{1}-\rho^4)(x_{1}-\rho^5)(x_{1}-\rho^6)\\
\\
x_{4}^{2}=(x_{1}-\rho)(x_{1}-\rho^3)(x_{1}-\rho^4)(x_{1}-\rho^5)\\
\\
z_{1}=\dfrac{1}{x_{1}}, \; z_{2}=\dfrac{x_{2}}{x_{1}^{2}}, \; z_{3}=\dfrac{\rho x_{3}x_{4}}{x_{1}^{2}(x_{1}-\rho^{4})(x_{1}-\rho^{5})}, \; z_{4}=\dfrac{\rho^{2} x_{3}}{x_{1}^{2}},\\
\\
w_{1}=x_{1}, \; w_{2}=x_{2}, \; w_{3}=x_{4}, \; w_{4}=\dfrac{x_{3}x_{4}}{(x_{1}-\rho^{4})(x_{1}-\rho^{5})},\\
\\
u_{1}=\dfrac{1}{x_{1}}, \; u_{2}=\dfrac{x_{2}}{x_{1}^{2}}, \; u_{3}=\dfrac{\rho^{2} x_{3}}{x_{1}^{2}}, \; u_{4}=\dfrac{\rho^{4} x_{4}}{x_{1}^{2}},\\
\\
v_{1}=x_{1}, \; v_{2}=x_{2}, \; v_{3}=\dfrac{x_{3}x_{4}}{(x_{1}-\rho^{4})(x_{1}-\rho^{5})}, \; v_{4}=x_{3},\\
\\
r_{1}=\dfrac{1}{x_{1}}, \; r_{2}=\dfrac{x_{2}}{x_{1}^{2}}, \; r_{3}=\dfrac{\rho^{4} x_{4}}{x_{1}^{2}}, \; r_{4}=\dfrac{\rho x_{3}x_{4}}{x_{1}^{2}(x_{1}-\rho^{4})(x_{1}-\rho^{5})}
\end{array}
\right\}
\end{equation}

\subsection{A permutation action}
Each $\tau \in \Gamma$ induces a natural bijection 
$$\widehat{\tau}:{{\mathbb C}}^{n} \to {{\mathbb C}}^{n}: (y_{1},\ldots, y_{n}) \mapsto (\sigma(y_{1}), \ldots, \sigma(y_{n})).$$

Let us consider the natural permutation action of $\Gamma$ on the coordinates of ${\mathbb C}^{24}$ defined by

$$\Theta(\sigma)(\vec{x},\vec{z},\vec{w},\vec{u},\vec{v},\vec{r})= (\vec{z},\vec{w},\vec{u},\vec{v},\vec{r},\vec{x})$$

Let us notice that the stabilizer $G$ of $\Phi(X)$ is trivial since
 $$G=\{\tau \in \Gamma: \Theta(\tau)(\Phi(X))=\Phi(X)\}=\{ \tau \in\Gamma: X^{\tau}=X\}=\{e\}.$$

Next, we should observe that, for each $\tau \in \Gamma$, the following diagram commutes \cite{HR}

\begin{equation}\label{diagrama1}
\begin{array}[c]{ccl} 
X&\stackrel{\Phi}{\rightarrow}&\Phi(X)\\ 
\downarrow\scriptstyle{f_{\tau}}&&\downarrow\scriptstyle{\Theta(\tau)}\\ 
X^{\tau}&\stackrel{\Phi^{\tau}}{\rightarrow}&\Theta(\tau)(\Phi(X))=\Phi^{\tau}(X^{\tau}) =\Phi(X)^{\tau}\\
\downarrow\scriptstyle{\widehat{\tau}^{-1}}&&\downarrow\scriptstyle{\widehat{\tau}^{-1}}\\ 
X&\stackrel{\Phi}{\rightarrow}&\Phi(X) 
\end{array}
\end{equation}

Similarly, it is not difficult to see that, for every $\eta,\tau \in\Gamma$, we have that 
$$(*) \quad \Theta(\eta) \circ \widehat{\tau}=\widehat{\tau} \circ \Theta(\eta).$$

\subsection{A birational isomorphism}
Assume we have computed a set of generators of the algebra of polynomials invariants under the linear action $\Theta(\Gamma)$, say $t_{1}$,..., $t_{N}$. 
It is not difficult to note that we may choose each $t_{j} \in {\mathbb Q}[\vec{x},\vec{z},\vec{w},\vec{u},\vec{v},\vec{r}]$. Next, we construct the rational map $$\Psi:{\mathbb C}^{24} \to {\mathbb C}^{N}$$
$$(\vec{x},\vec{z},\vec{w},\vec{u},\vec{v},\vec{r}) \mapsto (t_{1}, \ldots, t_{N})$$

\s

It can be checked, for each $\tau \in \Gamma$,  the following equalities:
\begin{equation}\label{eq4}
\begin{array}{l}
\Psi^{\tau}=\Psi\\
\Psi \circ \Theta(\tau)=\Psi
\end{array}
\end{equation}

Also, it holds (as we have chosen a set of generators of the invariant polynomials for the action of $\Theta(\Gamma)$) that $\Psi$ is a branched regular cover with Galois group $\Gamma$.

\s

It turns out that, if we set $W=\Psi(\Phi(X))$ and $L=\Psi \circ \Psi$, then 
$$L:X \to W$$ is a birational isomorphism (since the stabilizer $G$ of $\Phi(X)$ is trivial).

\subsection{$W$ can be defined by polynomials over ${\mathbb Q}$}
If $\tau \in \Gamma$, then
$$W^{\tau}=L(X)^{\tau}=L^{\tau}(X^{\tau})=\Psi^{\tau} \circ \Phi^{\tau} (X^{\tau})=\Psi \circ \Theta(\tau) (\Phi(X))=\Psi \circ \Phi(X)=L(X)=W.$$

This last set of equalities asserts that $W$ can be defined by polynomials with coefficient over ${\mathbb Q}$. In fact, it is almost clear that (using $L$ and the equations for $X$) that $W$ is defined by a set of polynomial equations over ${\mathbb Q}(\rho)$; say by $P_{1}, \ldots, P_{r} \in {\mathbb Q}(\rho)[t_{1},\ldots,t_{24}]$. The above set of equalities is telling us that, for each $\tau \in \Gamma$, the polynomials  
$P_{1}^{\tau}, \ldots, P_{r}^{\tau} \in {\mathbb Q}(\rho)[t_{1},\ldots,t_{24}]$ also are null over $W$. Let us consider the basis $\{1,\rho,\rho^{2},\rho^{3},\rho^{4},\rho^{5}\}$ of the ${\mathbb Q}$-vector space ${\mathbb Q}(\rho)$. The polynomials 
$${\rm Tr}(P_{j}), {\rm Tr}(\rho P_{j}), {\rm Tr}(\rho^{2} P_{j}), {\rm Tr}(\rho^{3} P_{j}), {\rm Tr}(\rho^{4} P_{j}), {\rm Tr}(\rho^{5} P_{j}) \in {\mathbb Q}[t_{1}, \ldots, t_{24}]$$
also define $W$ as their common zero locus, where 
$${\rm Tr}(aQ)=\sum_{l=0}^{5} \sigma^{l}(a) Q^{\sigma^{l}}$$

\subsection{In search of a set of invariant polynomials}
The program MAGMA \cite{MAGMA} may be used to obtain an explicit generating set of invariant polynomials. Unfortunately, in our case this was not possible in my computer (Mac OSX 10.5.8 2X2.66 GHz Dual-Core Intel Xeon 3GB 667 MHz DDR2 FB-DIMM).

By hand, we may find a collection of invariant polynomials given by

$$t_{j}=x_{j}+z_{j}+w_{j}+u_{j}+v_{j}+r_{j}; j=1,2,3,4$$
$$t_{4+j}=x_{j}^{2}+z_{j}^{2}+w_{j}^{2}+u_{j}^{2}+v_{j}^{2}+r_{j}^{2}; j=1,2,3,4$$
$$t_{8+j}=x_{j}^{3}+z_{j}^{3}+w_{j}^{3}+u_{j}^{3}+v_{j}^{3}+r_{j}^{3}; j=1,2,3,4$$
but I am not sure if they form a complete set of generators.

So, the idea will be to work in two steps. The first one will consist in produce a model for Fricke-Macbeath's curve in the degree two extension ${\mathbb Q}(\sqrt{-7})$ and then to use that model to obtain a model over ${\mathbb Q}$.  In the next two sections we describe these two steps.

\section{Step 1: A model of Fricke-Macbeath's curve over ${\mathbb Q}(\sqrt{-7})$} \label{Sec:Z2}
In this section we will explain how to obtain an explicit model for Fricke-Macbeath's curve over the quadratic extension ${\mathbb Q}(\sqrt{-7})$. Mostly of the computations have been carry out with MAGMA \cite{MAGMA} and with MATHEMATICA \cite{mathematica}.

Let $\tau=\sigma^{2}$ (an element of order $3$) and ${\mathbb Z}_{3} \cong N=\langle \tau \rangle \lhd \Gamma$. The fixed field of $N$ is ${\mathbb Q}(\rho+\rho^{2}+\rho^{4})={\mathbb Q}(\sqrt{-7})$.

Let us consider the rational map
$$\Phi_{1}:X \to {\mathbb C}^{12}$$
$$(x,y_{1},y_{2},y_{4}) \mapsto (\vec{x},\vec{w},\vec{v})$$

Again, we may see that $\Phi_{1}$ produces a birational isomorphism between $X$ and $\Phi_{1}(X)$ (its inverse is just the projection on the $\vec{x}$-coordinate; just as for $\Phi$).

Equations defining the algebraic curve $\Phi_{1}(X)$ are the following ones

\begin{equation}
\Phi_{1}(X)=\left\{\begin{array}{lcl}
x_{2}^{2}=(x_{1}-1)(x_{1}-\rho^3)(x_{1}-\rho^5)(x_{1}-\rho^6)\\
\\
x_{3}^{2}=(x_{1}-\rho^{2})(x_{1}-\rho^4)(x_{1}-\rho^5)(x_{1}-\rho^6)\\
\\
x_{4}^{2}=(x_{1}-\rho)(x_{1}-\rho^3)(x_{1}-\rho^4)(x_{1}-\rho^5)\\
\\
w_{1}=x_{1}, \; w_{2}=x_{2}, \; w_{3}=x_{4}, \; w_{4}=\dfrac{x_{3}x_{4}}{(x_{1}-\rho^{4})(x_{1}-\rho^{5})},\\
\\
v_{1}=x_{1}, \; v_{2}=x_{2}, \; v_{3}=\dfrac{x_{3}x_{4}}{(x_{1}-\rho^{4})(x_{1}-\rho^{5})}, \; v_{4}=x_{3}
\end{array}
\right\}
\end{equation}

\s

Next, we consider the linear permutation action of $N$ on the coordinates of ${\mathbb C}^{12}$ defined by
$$\Theta_{1}(\tau)(\vec{x},\vec{w},\vec{v})=(\vec{w},\vec{v},\vec{x})$$

Again, the stabilizer of $\Phi_{1}(X)$ is just the trivial group.

A generating set of invariant polynomials for the previous linear action can be obtained with MAGMA as

$$\begin{array}{lcl}
t_{1}&=&x_{1} + w_{1} + v_{1}\\
t_{2}&=&x_{2} + w_{2} + v_{2}\\
t_{3}&=&x_{3} + w_{3} + v_{3}\\
t_{4}&=&x_{4} + w_{4} + v_{4}\\
t_{5}&=&x_{1}^2 + w_{1}^2 + v_{1}^2\\
t_{6}&=&x_{2}^2 + w_{2}^2 + v_{2}^2\\
t_{7}&=&x_{3}^2 + w_{3}^2 + v_{3}^2\\
t_{8}&=&x_{4}^2 + w_{4}^2 + v_{4}^2\\
t_{9}&=&x_{1}^3 + w_{1}^3 + v_{1}^3\\
t_{10}&=&x_{2}^3 + w_{2}^3 + v_{2}^3\\
t_{11}&=&x_{3}^3 + w_{3}^3 + v_{3}^3\\
t_{12}&=&x_{4}^3 + w_{4}^3 + v_{4}^3
\end{array}
$$

The map
$$\Psi_{1}:{\mathbb C}^{12} \to {\mathbb C}^{12}$$
$$(\vec{x},\vec{w},\vec{v}) \mapsto (t_{1},...,t_{12})$$
satisfies the following properties: 
\begin{equation}
\begin{array}{l}
\Psi_{1}^{\tau^{j}}=\Psi_{1}, \; j=0,1,2;\\
\Psi_{1} \circ \Theta_{1}(\tau^{j})=\Psi_{1}, \; j=0,1,2.
\end{array}
\end{equation}

Also, it holds (as we have chosen a set of generators of the invariant polynomials for the action of $\Theta_{1}(N)$) that $\Psi_{1}$ is a branched regular cover with Galois group $N$. As in the previous case, it turns out that, if we set $Z_{1}=\Psi_{1}(\Phi_{1}(X))$ and $L_{1}=\Psi_{1} \circ \Psi_{1}$, then 
$$L_{1}:X \to Z_{1}$$ is a birational isomorphism (since the stabilizer of $\Phi_{1}(X)$ is trivial). 

Similarly, if $\eta \in N$, then
$$Z_{1}^{\eta}=L_{1}(X)^{\eta}=L_{1}^{\eta}(X^{\eta})=\Psi_{1}^{\eta} \circ \Phi_{1}^{\eta} (X^{\eta})=\Psi_{1} \circ \Theta_{1}(\eta) (\Phi_{1}(X))=\Psi_{1} \circ \Phi_{1}(X)=L_{1}(X)=Z,$$
so $Z_{1}$ can be defined by polynomials with coefficient over ${\mathbb Q}(\sqrt{-7})$.  

We have tried to compute directly equations for $Z_{1}$ using MAGMA and the explicit form of $L_{1}$, but unfortunately we couldn't get an answer (it is very heavy computational task). We will proceed in the search of such equation by hands.

It is clear that
$$x_{1}=\frac{t_{1}}{3}$$
$$x_{2}=\frac{t_{2}}{3}$$
$$t_{4}=t_{3}$$
$$(*) \; x_{4}=\frac{(t_{3}-x_{3})(\frac{t_{1}}{3}-\rho^{4})(\frac{t_{1}}{3}-\rho^{5})}{x_{3}+(\frac{t_{1}}{3}-\rho^{4})(\frac{t_{1}}{3}-\rho^{5})}$$
$$t_{5}=\frac{t_{1}^{2}}{3}$$
$$t_{6}=\frac{t_{2}^{2}}{3}$$
$$t_{8}=t_{7}$$
$$(**) \; x_{4}^{2}=\frac{(t_{7}-x_{3}^{2})(\frac{t_{1}}{3}-\rho^{4})^{2}(\frac{t_{1}}{3}-\rho^{5})^{2}}{x_{3}^{2}+(\frac{t_{1}}{3}-\rho^{4})^{2}(\frac{t_{1}}{3}-\rho^{5})^{2}}$$
$$t_{9}=\frac{t_{1}^{3}}{9}$$
$$t_{10}=\frac{t_{2}^{3}}{9}$$
$$t_{12}=t_{11}$$
$$(***)\; x_{4}^{3}=\frac{(t_{11}-x_{3}^{3})(\frac{t_{1}}{3}-\rho^{4})^{3}(\frac{t_{1}}{3}-\rho^{5})^{3}}{x_{3}^{3}+(\frac{t_{1}}{3}-\rho^{4})^{3}(\frac{t_{1}}{3}-\rho^{5})^{3}}$$

Equality $(*)$ permits to obtain $x_{4}$ uniquely in terms of $t_{1}$ and $x_{3}$.

Equation $$x_{2}^{2}=(x_{1}-1)(x_{1}-\rho^{3})(x_{1}-\rho^{5})(x_{1}-\rho^{6})$$ asserts the polynomial equation (relating $t_{1}$ and $t_{2}$)
$$P_{1}(t_{1},t_{2},t_{3},t_{7},t_{11})=
-81+27(1+(\rho+\rho^2+\rho^4))t_{1}+9t_{1}^{2}-3(\rho+\rho^2+\rho^4)t_{1}^{3}-t_{1}^{4}+9t_{2}^{2}=0.$$

Notice that $P_{1}(t_{1},t_{2},t_{3},t_{7},t_{11}) \in {\mathbb Q}(\sqrt{-7})[t_{1},t_{2},t_{3},t_{7},t_{11}]$.

From the equation $$x_{3}^{2}=(x_{1}-\rho^{2})(x_{1}-\rho^{4})(x_{1}-\rho^{5})(x_{1}-\rho^{6})$$ we obtain
the equation
$$(1) \; x_{3}^{2}=(t_{1}-3\rho^{2})(t_{1}-3\rho^{4})(t_{1}-3\rho^{5})(t_{1}-3\rho^{6})/81$$

From the equation $$x_{4}^{2}=(x_{1}-\rho)(x_{1}-\rho^{3})(x_{1}-\rho^{4})(x_{1}-\rho^{5})$$ we obtain
the equation
$$(2) \; x_{4}^{2}=(t_{1}-3\rho)(t_{1}-3\rho^{3})(t_{1}-3\rho^{4})(t_{1}-3\rho^{5})/81$$

In this way, by replacing the above values for $x_{3}^{2}$ and $x_{4}^{2}$ (obtained in $(1)$ and $(2)$) in the above equality $(**)$, we obtain the polynomial equation 
$$P_{2}(t_{1},t_{2},t_{3},t_{7},t_{11})=
27+27 (\rho+ \rho^2+ \rho^4)-18 t_{1}-3(1 + (\rho + \rho^2 + \rho^4)) t_{1}^2-2 t_{1}^3-t_{1}^4+27 t_{7}=0.
$$

Notice that $P_{2}(t_{1},t_{2},t_{3},t_{7},t_{11}) \in {\mathbb Q}(\sqrt{-7})[t_{1},t_{2},t_{3},t_{7},t_{11}]$.

Also, if we replace, in equality $(***)$ the values of $x_{3}^{3}$ by $x_{3}(x_{1}-\rho^{2})(x_{1}-\rho^{4})(x_{1}-\rho^{5})(x_{1}-\rho^{6})/81$ and $x_{4}^{3}$ by $x_{4}(t_{1}-3\rho)(t_{1}-3\rho^{3})(t_{1}-3\rho^{4})(t_{1}-3\rho^{5})/81$, where $x_{4}$ is given in $(*)$, then we obtain a polynomial which is of degree one in the variable $x_{3}$.

\s

$x_{3}=(-9 \rho^2 (-162 t_{1} -18 t_{1}^3 + 4 t_{1}^5-243 (1+t_{11}) + t_{1}^2 (27 - 54 t_{3}) + 6 t_{1}^4 t_{3}) + 3 (729 + 18 t_{1}^4 - 6 t_{1}^5 - 27 t_{1}^3 (-6 + t_{3}) - t_{1}^6 (-2 + t_{3}) + 243 t_{1} (3 + t_{3}) + 81 t_{1}^2 (2 + t_{11} + t_{3}))+\rho^3 (2187 - t_{1}^7 + 27 t_{1}^4 (-6 + t_{3}) + 9 t_{1}^5 (-3 + t_{3}) + 486 t_{1}^2 t_{3} + 81 t_{1}^3 (1 + t_{3}) + 729 t_{1} (1 + 2 t_{3})) + \rho^5 (2187 + 27 t_{1}^4 + 12 t_{1}^6 + t_{1}^7 - 729 t_{1} (-1 + t_{11} - t_{3}) + 729 t_{1}^2 t_{3} + 81 t_{1}^3 (5 + t_{3}) + 9 t_{1}^5 (1 + 2 t_{3})) + \rho (2916 t_{1} + 3 t_{1}^6 - t_{1}^7 - 81 t_{1}^3 (-6 + t_{3}) - 2187 (-2 + t_{3}) - 27 t_{1}^4 (-2 + t_{3}) + 9 t_{1}^5 (2 + t_{3}) + 243 t_{1}^2 (5 + 2 t_{3})) + \rho^4 (2187 + t_{1}^7 - 729 t_{1} (-3 + t_{11} - 2 t_{3}) - 81 t_{1}^3 (-1 + t_{3}) + 27 t_{1}^4 (1 + t_{3}) + 9 t_{1}^5 (-1 + 2 t_{3}) + 243 t_{1}^2 (1 + 3 t_{3})))/(9 (t_{1}^5 - 243 t_{11} + 27 t_{1}^2 (-1 + t_{3}) + 81 t_{1} t_{3} + 9 t_{1}^3 t_{3} + 3 t_{1}^4 t_{3} + \rho (3 + t_{1}) (-81 + 18 t_{1}^2-9 t_{1}^3 + 2 t_{1}^4 + 27 t_{1} t_{3}) + 27 \rho^2 t_{1} (3 + t_{1}^2 + t_{1} (3 + t_{3})) + \rho^4 t_{1} (243 + 3 t_{1}^3 + t_{1}^4 + 9 t_{1}^2 (-1 + t_{3}) + 27 t_{1} (3 + t_{3})) + \rho^5 (-6 t_{1}^4 + t_{1}^5 + 243 (1 + t_{3}) + 81 t_{1} (2 + t_{3}) + 9 t_{1}^3 (2 + t_{3}) + 27 t_{1}^2 (3 + t_{3})) + \rho^3 t_{1} (162 + 36 t_{1}^2 + 6 t_{1}^3 + 2 t_{1}^4 + 27 t_{1} (4 + t_{3}))))
$

\s

Then, using $(*)$, we obtain

\s

$x_{4}=-((3 \rho^4-t_{1}) (3 \rho^5-t_{1}) (-\rho^3 (-2187-729 t_{1} + t_{1}^7 + 243 t_{1}^2 t_{3} (2 + t_{3}) + 9 t_{1}^5 (3 + t_{3}) + 27 t_{1}^4 (6 + t_{3}) + 81 t_{1}^3 (-1 + 3 t_{3})) + \rho^4 (2187 + 27 t_{1}^4 + t_{1}^7 + 9 t_{1}^5 (-1 + t_{3})-729 t_{1} (-3 + t_{11} + t_{3})-243 t_{1}^2 (-1 + t_{3}^2)-81 t_{1}^3 (-1 + t_{3}^2)) + \rho (4374 + 486 t_{1}^3 + 54 t_{1}^4 + 3 t_{1}^6-t_{1}^7-9 t_{1}^5 (-2 + t_{3})-243 t_{1}^2 (-5 + t_{3}^2)-729 t_{1} (-4-t_{3} + t_{3}^2))-3 (t_{1}^6 (-2 + t_{3}) + 3 t_{1}^5 (2 + t_{3})-729 (1 + t_{11} t_{3})-81 t_{1}^2 (2 + t_{11} + 2 t_{3}-t_{3}^2) + 9 t_{1}^4 (-2 + t_{3}^2) + 243 t_{1} (-3-t_{3} + t_{3}^2) + 27 t_{1}^3 (-6 + t_{3} + t_{3}^2))-9 \rho^2 (4 t_{1}^5-243 (1 + t_{11}) + 81 t_{1} (-2 + t_{3}) + 6 t_{1}^4 t_{3} + 9 t_{1}^3 (-2 + 3 t_{3}) + 27 t_{1}^2 (1 + t_{3} + t_{3}^2)) + \rho^5 (12 t_{1}^6 + t_{1}^7-243 t_{1}^2 t_{3}^2 + 9 t_{1}^5 (1 + t_{3}) + 27 t_{1}^4 (1 + 2 t_{3})-81 t_{1}^3 (-5 + t_{3} + t_{3}^2)-2187 (-1 + t_{3} + t_{3}^2)-729 t_{1} (-1 + t_{11} + t_{3} + t_{3}^2))))/(9 (567 t_{1}^3 + 6 t_{1}^6 + t_{1}^7 + \rho (-3 + t_{1}) (-54 t_{1}^3 + t_{1}^6 + 9 t_{1}^4 (-4 + t_{3}) + 729 (-2 + t_{3}) + 243 t_{1} (-2 + t_{3})-81 t_{1}^2 (-2 + t_{3})) + 27 t_{1}^4 (-7 + t_{3}) + 9 t_{1}^5 (-5 + t_{3}) + 2187 (2 + t_{3}) + 243 t_{1}^2 (-1 + 2 t_{3}) + 729 t_{1} (1 + 2 t_{3}) + \rho^5 (2187 + 216 t_{1}^4 + 3 t_{1}^6 + 2 t_{1}^7 + 729 t_{1} t_{3} + 729 t_{1}^2 (1 + t_{3}) + 18 t_{1}^5 (2 + t_{3}) + 81 t_{1}^3 (16 + t_{3})) + \rho^3 t_{1} (9 t_{1}^5 + t_{1}^6 + 27 t_{1}^3 (-4 + t_{3}) + 9 t_{1}^4 (3 + t_{3}) + 81 t_{1}^2 (5 + t_{3}) + 729 (-5 + 2 t_{3}) + 243 t_{1} (-3 + 2 t_{3})) + \rho^4 (2187 + 6 t_{1}^6 + 2 t_{1}^7-81 t_{1}^3 (-14 + t_{3}) + 18 t_{1}^5 (-2 + t_{3}) + 1458 t_{1} t_{3} + 27 t_{1}^4 (5 + t_{3}) + 243 t_{1}^2 (1 + 3 t_{3}))-9 \rho^2 (-243 + 243 t_{1}-27 t_{1}^3 + t_{1}^5-54 t_{1}^2 (-5 + t_{3}) + t_{1}^4 (-9 + 6 t_{3}))))$.

\s

Now, using such values for $x_{3}$ and $x_{4}$, and replacing them in $(1)$ and $(2)$ above, we obtain another two polynomials identities $P_{1}(t_{1},t_{3},t_{7},t_{11})=0$ and $P_{2}(t_{1},t_{3},t_{7},t_{11})=0$, where the polynomials are defined over ${\mathbb Q}(\rho)$ (see the Appendix). In this way,

$$Z_{1}=\left\{ \begin{array}{c}
t_{4}=t_{3}, \; 
3t_{5}=t_{1}^{2}, \;
3t_{6}=t_{2}^{2}, \;
t_{8}=t_{9}\\
9t_{9}=t_{1}^{3}, \;
9t_{10}=t_{2}^{3}, \;
t_{12}=t_{11}\\
P_{1}(t_{1},t_{2},t_{3},t_{7},t_{11})=0\\
P_{2}(t_{1},t_{2},t_{3},t_{7},t_{11})=0\\
P_{3}(t_{1},t_{2},t_{3},t_{7},t_{11})=0\\
P_{4}(t_{1},t_{2},t_{3},t_{7},t_{11})=0
\end{array}
\right\} \subset {\mathbb C}^{12}
$$

Notice that, by the above computations, we have explicitly the inverse of $L_{1}$ given as
$$L_{1}^{-1}:Z_{1} \to X$$
$$(t_{1},...,t_{12}) \mapsto (x_{1},x_{2},x_{3},x_{4}),$$
that is, we may write $x_{1}$, $x_{2}$, $x_{3}$ and $x_{4}$ in terms of $t_{1}$,..., $t_{12}$ (in fact, only in terms of $t_{1}$, $t_{2}$, $t_{3}$, $t_{7}$ and $t_{11}$). 

As the variables $t_{1}$,..., $t_{12}$ are uniquely determined only by the variables $t_{1}$, $t_{2}$, $t_{3}$, $t_{7}$ and $t_{11}$, if we consider the projection
$$\pi:{\mathbb C}^{12} \to {\mathbb C}^{5}$$ 
$$(t_{1},...,t_{12}) \mapsto (t_{1},t_{2},t_{3},t_{7},t_{11}),$$
then 

$$L_{1}^{*}=\pi \circ L_{1}:X \to Z_{2}$$
$$L_{1}^{*}(x,y_{1},y_{2},y_{4})$$
$$||$$
$$
\left(3x,3y_{1},y_{2}+y_{4}+\frac{y_{2}y_{4}}{(x-\rho^{4})(x-\rho^{5})}, y_{2}^{2}+y_{4}^{2}+\frac{y_{2}^{2}y_{4}^{2}}{(x-\rho^{4})^{2}(x-\rho^{5})^{2}},  y_{2}^{3}+y_{4}^{3}+\frac{y_{2}^{3}y_{4}^{3}}{(x-\rho^{4})^{3}(x-\rho^{5})^{3}} \right)
$$
is a birational isomorphism, where 
$$Z_{2}=\left\{ \begin{array}{c}
P_{1}(t_{1},t_{2},t_{3},t_{7},t_{11})=0\\
P_{2}(t_{1},t_{2},t_{3},t_{7},t_{11})=0\\
P_{3}(t_{1},t_{2},t_{3},t_{7},t_{11})=0\\
P_{4}(t_{1},t_{2},t_{3},t_{7},t_{11})=0
\end{array}
\right\} \subset {\mathbb C}^{5}
$$ 

The inverse $(L_{1}^{*})^{-1}:Z_{2} \to X$ is given as
$$(L_{1}^{*})^{-1}(t_{1},t_{2},t_{3},t_{7},t_{11})=(x_{1},x_{2},x_{3},x_{4}),$$
where $x_{1}$,..., $x_{4}$ are given by the previous formulae.

As $Z_{1}^{\eta}=Z_{1}$, for every $\eta \in N$, as as $\pi$ is defined over ${\mathbb Q}$, we see that 
$Z_{2}^{\eta}=Z_{2}$, for every $\eta \in N$, that is, $Z_{2}$ can be defined by polynomials over ${\mathbb Q}(\sqrt{-7})$.  To obtain such equations over ${\mathbb Q}(\sqrt{-7})$, we replace each polynomial $P_{j}$ ($j=3,4$) by the polynomials (with coefficients in ${\mathbb Q}(\sqrt{-7})$)
$$Q_{j,1}={\rm Tr}(P_{j}), \; Q_{j,2}={\rm Tr}(\rho P_{j}), \; Q_{j,3}={\rm Tr}(\rho^{2} P_{j})$$
that is

$$Z_{2}=\left\{ \begin{array}{c}
P_{1}(t_{1},t_{2},t_{3},t_{7},t_{11})=0\\
P_{2}(t_{1},t_{2},t_{3},t_{7},t_{11})=0\\
P_{3}(t_{1},t_{2},t_{3},t_{7},t_{11})+P_{3}(t_{1},t_{2},t_{3},t_{7},t_{11})^{\tau}+P_{3}(t_{1},t_{2},t_{3},t_{7},t_{11})^{\tau^{2}}=0\\
\rho P_{3}(t_{1},t_{2},t_{3},t_{7},t_{11})+\rho^{2} P_{3}(t_{1},t_{2},t_{3},t_{7},t_{11})^{\tau}+\rho^{4} P_{3}(t_{1},t_{2},t_{3},t_{7},t_{11})^{\tau^{2}}=0\\
\rho^{2} P_{3}(t_{1},t_{2},t_{3},t_{7},t_{11})+\rho^{4} P_{3}(t_{1},t_{2},t_{3},t_{7},t_{11})^{\tau}+\rho P_{3}(t_{1},t_{2},t_{3},t_{7},t_{11})^{\tau^{2}}=0\\
P_{4}(t_{1},t_{2},t_{3},t_{7},t_{11})+P_{4}(t_{1},t_{2},t_{3},t_{7},t_{11})^{\tau}+P_{4}(t_{1},t_{2},t_{3},t_{7},t_{11})^{\tau^{2}}=0\\
\rho P_{4}(t_{1},t_{2},t_{3},t_{7},t_{11})+\rho^{2} P_{4}(t_{1},t_{2},t_{3},t_{7},t_{11})^{\tau}+\rho^{4} P_{4}(t_{1},t_{2},t_{3},t_{7},t_{11})^{\tau^{2}}=0\\
\rho^{2} P_{4}(t_{1},t_{2},t_{3},t_{7},t_{11})+\rho^{4} P_{4}(t_{1},t_{2},t_{3},t_{7},t_{11})^{\tau}+\rho P_{4}(t_{1},t_{2},t_{3},t_{7},t_{11})^{\tau^{2}}=0
\end{array}
\right\} \subset {\mathbb C}^{5}
$$

We have obtained an explicit model $Z_{2}$ for Fricke-Macbeath's curve over ${\mathbb Q}(\sqrt{-7})$ together explicit isomorphisms $L_{1}^{*}:X \to Z_{2}$ and $(L_{1}^{*})^{-1}:Z_{2} \to X$.

\s
\noindent
\begin{rema}\label{dibujo}
Notice that the regular dessin d'enfants $(X,\beta)$, given in the Introduction, is isomorphic to that provided by the pair $(Z_{2},\beta^{*})$, where $\beta^{*}(t_{1},t_{2},t_{3},t_{7},t_{11})=(t_{1}/3)^{7}$; that is, such a dessin d'enfants is defined over ${\mathbb Q}(\sqrt{-7})$. 
\end{rema}

\s

\section{Step 2: A model of Fricke-Macbeath's curve over ${\mathbb Q}$}
Let us now consider the explicit model $Z_{2} \subset {\mathbb C}^{5}$ over ${\mathbb Q}(\sqrt{-7})$, constructed in the previous section. We keep the same notations.
Let $M={\rm Gal}({\mathbb Q}(\sqrt{-7})/{\mathbb Q})=\langle \eta \rangle \cong {\mathbb Z}_{2}$, where $\eta$ is the complex conjugation.

As already noticed in the Introduction,  $X$ admits the following anticonformal involution
$$J(x,y_{1},y_{2},y_{4})=\left(\dfrac{1}{\overline{x}}, \dfrac{\overline{y_{1}}}{\overline{x}^{2}},
\dfrac{\rho^{5} \overline{y_{2}}}{\overline{x}^{2}},\dfrac{\rho^{3} \overline{y_{4}}}{\overline{x}^{2}}\right)$$

In this way, $T=L_{1}^{*} \circ J \circ (L_{1}^{*})^{-1}$ is an anticonformal involution of $Z_{2}$. It is not difficult to see that by setting $g_{e}=I$ and $g_{\eta}=S \circ T$, where $S(t_{1},t_{2},t_{3},t_{7},t_{11})=(\overline{t_{1}}, \overline{t_{2}},\overline{t_{3}},\overline{t_{7}},\overline{t_{11}})$, we obtain a Weil's datum for $Z_{2}$.

Now, identically as done before, we consider the rational map
$$\Phi_{2}:Z_{2} \to {\mathbb C}^{10}$$
$$(t_{1},t_{2},t_{3},t_{7},t_{11}) \mapsto (t_{1},t_{2},t_{3},t_{7},t_{11},s_{1},s_{2},s_{3},s_{7},s_{11})$$
where $g_{\eta}(t_{1},t_{2},t_{3},t_{7},t_{11})=(s_{1},s_{2},s_{3},s_{7},s_{11})$. Then, again we may see that
$\Phi_{2}$ induces a birational isomorphism between $Z_{2}$ and $\Phi_{2}(Z_{2})$.

In this case, 
$$\Phi_{2}(Z_{2})=\left\{\begin{array}{c}
Q_{1,1}(t_{1},t_{2},t_{3},t_{7},t_{11})=\cdots= Q_{4,3}(t_{1},t_{2},t_{3},t_{7},t_{11})=0\\
g_{\eta}(t_{1},t_{2},t_{3},t_{7},t_{11})=(s_{1},s_{2},s_{3},s_{7},s_{11})
\end{array}
\right\} \subset {\mathbb C}^{10}.$$

In this case, the Galois group $M$ induces the permutation action $\Theta_{2}(M)$ defined as
$$\Theta(\eta)(t_{1},t_{2},t_{3},t_{7},t_{11},s_{1},s_{2},s_{3},s_{7},s_{11})=(s_{1},s_{2},s_{3},s_{7},s_{11},t_{1},t_{2},t_{3},t_{7},t_{11})$$

A set of generators for the invariant polynomials (with respect to the previous permutation action) is given by

$$q_{1}=t_{1}+s_{1}, \; q_{2}=t_{2}+s_{2}, \; q_{3}=t_{3}+s_{3},$$
$$q_{4}=t_{7}+s_{7},\; q_{5}=t_{11}+s_{11},\; q_{6}=t_{1}^{2}+s_{1}^{2},$$
$$q_{7}=t_{2}^{2}+s_{2}^{2},\; q_{8}=t_{3}^{2}+s_{3}^{2}, \;q_{9}=t_{7}^{2}+s_{7}^{2},$$
$$q_{10}=t_{11}^{2}+s_{11}^{2}$$

Then the rational map
$$\Psi_{2}:{\mathbb C}^{10} \to {\mathbb C}^{10}$$
$$(t_{1},t_{2},t_{3},t_{7},t_{11},s_{1},s_{2},s_{3},s_{7},s_{11}) \mapsto (q_{1},q_{2},q_{3},q_{4},q_{5},q_{6},q_{7}q_{8},q_{9},q_{10})$$
satisfies the following properties: 
\begin{equation}
\begin{array}{l}
\Psi_{2}^{\eta}=\Psi_{2};\\
\Psi_{2} \circ \Theta_{2}(\eta)=\Psi_{2}.
\end{array}
\end{equation}

There are two possibilities: 
\begin{enumerate}
\item $\Phi_{2}(Z_{2})=\Phi_{2}(Z_{2})$; in which case $Z_{2}^{\eta}=Z_{2}$ and $Z_{2}$ will be already defined over ${\mathbb Q}$ (which seems not to be the case); and
\item the stabilizer of $\Phi_{2}(Z_{2})$ under $\Theta_{2}(M)$ is trivial.
\end{enumerate}

Under the assumption (2) above, we have that $\Psi_{2}:\Phi_{2}(Z_{2}) \to W=\Psi_{2}(\Phi_{2}(Z_{2}))$ is a birational isomorphism and that, as before, $W$ is defined over ${\mathbb Q}$. That is, the map
$L_{2}=\Psi_{2} \circ \Phi_{2}:Z_{2} \to W$ is an explicit birational isomorphism and $W$ is defined over ${\mathbb Q}$. In this way, $L=L_{2} \circ L_{1}^{*}:X \to W$ is an explicit birational isomorphism as desired at the beginning.

\s

Now, again as in the previous section, as $R_{2}$ and $Z_{2}$ are explicitly given, one may compute explicit equations for $W$ over ${\mathbb Q}(\sqrt{-7})$, say by the polynomials $A_{1},...,A_{m} \in {\mathbb Q}(\sqrt{-7})[q_{1},...,q_{10}]$ (this may be done with MAGMA \cite{MAGMA} or by hands, but computations are heavy and very long). We do not write these equations as they are large expressions to be written down in this paper (anyway, all of this can be done explicitly with the computer). To obtain equations over ${\mathbb Q}$ we replace each $A_{j}$ (which is not already defined over ${\mathbb Q}$) by the traces $A_{j}+A_{j}^{\eta}$, $iA_{j}-iA_{j}^{\eta}$. 

\section{A remark on the elliptic curves in \eqref{eq0}}\label{appendix}

\subsection{A connection to homology covers}
Let us set $\lambda_{1}=1$, $\lambda_{2}=\rho$, $\lambda_{3}=\rho^{2}$, $\lambda_{4}=\rho^{3}$, $\lambda_{5}=\rho^{4}$, $\lambda_{6}=\rho^{5}$ and $\lambda_{7}=\rho^{6}$.
It is known that if $S$ is Fricke-Macbeath's curve, then it admits a regular branched cover $Q:S \to \widehat{\mathbb C}$ whose branch locus is the set $\{\lambda_{1},\lambda_{2},\lambda_{3},\lambda_{4},\lambda_{5},\lambda_{6},\lambda_{7}\}$ and deck group $G \cong {\mathbb Z}_{2}^{3}$. Let us consider a Fuchsian group 
$$\Gamma=\langle \alpha_{1},...,\alpha_{7}: \alpha_{1}^{2}=\cdots=\alpha_{7}^{2}=\alpha_{1} \alpha_{2} \cdots \alpha_{7}=1\rangle$$
acting on the hyperbolic plane ${\mathbb H}^{2}$ uniformizing the orbifold $S/G$. 

If $\Gamma'$ denotes the derived subgroup of $\Gamma$, then $\Gamma'$ acts freely and $\widehat{S}={\mathbb H}^{2}/\Gamma'$ is a closed Riemann surface. Let $H=\Gamma/\Gamma' \cong {\mathbb Z}_{2}^{6}$; a group of conformal automorphisms of $\widehat{S}$.
Then there exists a set of generators of $H$, say $a_{1}$,..., $a_{6}$, so that the only elements of $H$ acting with fixed points are these and $a_{7}=a_{1}a_{2}a_{3}a_{4}a_{5}a_{6}$. In fact, $\widehat{S}$ corresponds to the algebraic curve (called the homology cover of $S/H$)
$$\widehat{S}=\left\{ \begin{array}{rcl}
x_{1}^{2}+x_{2}^{2}+x_{3}^{2}&=&0\\
\left(\dfrac{\lambda_{3}-1}{\lambda_{4}-1}\right) x_{1}^{2}+x_{2}^{2}+x_{4}^{2}&=&0\\
\\
\left(\dfrac{\lambda_{4}-1}{\lambda_{5}-1}\right) x_{1}^{2}+x_{2}^{2}+x_{5}^{2}&=&0\\
\\
\left(\dfrac{\lambda_{5}-1}{\lambda_{6}-1}\right) x_{1}^{2}+x_{2}^{2}+x_{6}^{2}&=&0\\
\\
\left(\dfrac{\lambda_{6}-1}{\lambda_{7}-1}\right) x_{1}^{2}+x_{2}^{2}+x_{7}^{2}&=&0
\end{array}
\right\} \subset {\mathbb P}_{\mathbb C}^{6},
$$
and $a_{j}$ is just multiplication by $-1$ at the coordinate $x_{j}$. The regular branched cover $P:\widehat{S} \to \widehat{\mathbb C}$ is given bt
$$P([x_{1}: x_{2}: x_{3}: x_{4}: x_{5}: x_{6}: : x_{7}])=\frac{x_{2}^{2}+x_{1}^{2}}{x_{2}^{2}+\lambda_{7} x_{1}^{2}}=z.$$

The branch locus of $P$ is given by the set of the $7th$-roots of unity $\{\lambda_{1},..., \lambda_{7}\}$. All details can be found in \cite{CGHR}.

By classical covering theory, there should be a subgroup $K<H$, $K \cong {\mathbb Z}_{2}^{3}$, acting freely on $\widehat{S}$ so that there is an isomorphism
$\phi:S \to \widehat{S}/K$ with $\phi G \phi^{-1}=H/K$.

Let us also observe that the rotation $R(z)=\rho z$ lifts under $P$ to an automorphism
$T$ of $\widehat{S}$ of order $7$ of the form
$$T([x_{1}:\cdots:x_{7}])=[c_{1}x_{7}:c_{2}x_{1}:c_{3}x_{2}:c_{4}x_{3}:c_{5}x_{4}:c_{6}x_{5}:c_{7}x_{6}]$$
for suitable comples numbers $c_{j}$. One has that $T a_{j} T^{-1}=a_{j+1}$, for $j=1,...,6$ and $T a_{7} T^{-1}=a_{1}$. The subgroup $K$ above must satisfies that $T K T^{-1}=K$ as $R$ lifts to Fricke-Macbeath's curve (as noticed in the Introduction).

\subsection{About the elliptic curves in Fricke-Macbeath's curve}
The subgroup $K^{*}=\langle a_{1}a_{3}a_{7}, a_{2}a_{3}a_{5},a_{1}a_{2}a_{4}\rangle \cong {\mathbb Z}_{2}^{3}$ acts freely on $\widehat{S}$ and it is normalized by  $T$. In particular, 
$S^{*}=\widehat{S}/K$ is a closed Riemann surface of genus $7$ admitting the group 
$L=H/K^{*}=\{e,a_{1}^{*},...,a_{7}^{*}\} \cong {\mathbb Z}_{2}^{3}$ (where $a_{j}^{*}$ is the involution induced by $a_{j}$) as a group of automorphisms and it also has an automorphism $T^{*}$ of order $7$ (induced by $T$) permuting cyclically the involutions $a_{j}^{*}$. As $S^{*}/\langle L, T^{*}\rangle=\widehat{S}/\langle H,T\rangle$ has signature $(0;2,7,7)$, we may see that $S=S^{*}$ and $K=K^{*}$.

We may see that $L=\langle a_{1}^{*},a_{2}^{*},a_{3}^{*}\rangle$ and $a_{4}^{*}=a_{1}^{*}a_{2}^{*}$, $a_{5}^{*}=a_{2}^{*}a_{3}^{*}$, $a_{6}^{*}=a_{1}^{*}a_{2}^{*}a_{3}^{*}$ and $a_{7}^{*}=a_{1}^{*}a_{3}^{*}$. In this way, we may see that every involution of $H/K$ is induced by one of the involutions (and only one) with fixed points; so every involutions in $L$ acts with $4$ fixed points.

 Let $a_{i}^{*}, a_{j}^{*} \in H/K$ be any two different involutions, so $\langle a_{i}^{*}, a_{j}^{*}  \rangle \cong {\mathbb Z}_{2}^{2}$. Then, by Riemann-Hurwitz, the quotient surface $S/\langle a_{i}^{*}, a_{j}^{*}  \rangle$ is a closed Riemann surface of genus $1$ with six cone points of order $2$. These six cone points are projected onto three of the cone points of $S/H$. These points are $\lambda_{i}$, $\lambda_{j}$ and $\lambda_{r}$, where $a_{r}^{*}=a_{i}^{*} a_{j}^{*}$. In this way, 
 the genus one surface is given by the elliptic curve 
 $$y^{2}=\prod_{k \notin \{i,j,r\}} (x-\lambda_{k})$$
 
 So, for instance, if we consider $i=2$ and $j=3$, then $r=5$ and the elliptic curve is
 $$y_{1}^{2}=(x-1)(x-\rho^3)(x-\rho^5)(x-\rho^6).$$
 
 If $i=1$ and $j=2$, then $r=4$ and the elliptic curve is 
$$y_{2}^{2}=(x-\rho^{2})(x-\rho^4)(x-\rho^5)(x-\rho^6).$$

If $i=1$ and $j=3$, then $r=7$ and the elliptic curve is
$$y_{4}^{2}=(x-\rho)(x-\rho^3)(x-\rho^4)(x-\rho^5).$$

We have obtained the three elliptic curves appearing in the Fricke-Macbeath's curve \eqref{eq0}.

 \subsection{Another model for Fricke-Macbeath's curve}
 The above description of Fricke-Macbeath's curve in terms of the homology cover $\widehat{S}$ permits to obtain an explicit model.
Let us consider now an affine model of $\widehat S$, say by taking $x_{7}=1$, with we denote by $\widehat{S}^{0}$. In this case the involution $a_{7}$ is multiplication of all coordinates by $-1$. A set of generators for the algebra of invariant polynomials in ${\mathbb C}[x_1,x_2,x_3,x_4,x_5,x_6]$ under the natural linear action induced by $K$ is
 $$t_1=x_{1}^{2}, t_2=x_{2}^{2}, t_3=x_{3}^{2}, t_4=x_{4}^{2}, t_5=x_{5}^{2}, t_6=x_{6}^{2}, t_7=x_{1}x_{2}x_{5}, t_8=x_{1}x_{2},x_{3}x_{6}$$
 $$t_{9}=x_{1}x_{4}x_{6}, t_{10}=x_{1}x_{3}x_{4}x_{5}, t_{11}=x_{2}x_{4}x_{5}x_{6}, t_{12}=x_{2}x_{3}x_{4}, t_{13}=x_{3}x_{5}x_{6}.$$
 
 If we set $$F:\widehat{S}^{0} \to {\mathbb C}^{13}$$
 $$F(x_1,x_2,x_3,x_4,x_5,x_6)=(t_{1}, t_{2}, t_{3}, t_{4}, t_{5}, t_{6}, t_{7}, t_{8},t_{9},t_{10},t_{11},t_{12},t_{13}),$$
 then $F(\widehat{S}^{0})$ will provide a model for  Fricke-Macbeath's curve $S$ (affine). Equations for such an affine model of $S$ are
 $$\left\{ \begin{array}{c}
 t_1+t_2+t_3=0\\
 \left(\frac{\lambda_3 - 1}{\lambda_4 -1}\right) t_1+t_2+t_4=0\\
  \left(\frac{\lambda_4 - 1}{\lambda_5 -1}\right) t_1+t_2+t_5=0\\
  \left(\frac{\lambda_5 - 1}{\lambda_6 -1}\right) t_1+t_2+t_6=0\\
  \left(\frac{\lambda_6 - 1}{\lambda_7 -1}\right) t_1+t_2+1=0\\
t_{6}t_{10} = t_{9}t_{13}, \; 
t_{6}t_{7}t_{12} = t_{8}t_{11},\;
t_{5}t_{9}t_{12} = t_{10}t_{11}\\
t_{5}t_{8} = t_{7}t_{13}, \;
t_{5}t_{6}t_{12} = t_{11}t_{13}, \;
t_{4}t_{8} = t_{9}t_{12} \\
t_{4}t_{7}t_{13} = t_{10}t_{11}, \;
t_{4}t_{6}t_{7} = t_{9}t_{11}, \;
t_{3}t_{11} = t_{12}t_{13} \\
t_{3}t_{6}t_{7} = t_{8}t_{13}, \;
t_{3}t_{5}t_{9} = t_{10}t_{13}, \;
t_{3}t_{5}t_{6} = t_{13}^{2}\\
t_{3}t_{4}t_{7} = t_{10}t_{12}, \;
t_{2}t_{10} = t_{7}t_{12}, \;
t_{2}t_{9}t_{13} = t_{8}t_{11} \\
t_{2}t_{5}t_{9} = t_{7}t_{11}, \;
t_{2}t_{4}t_{13} = t_{11}t_{12}, \;
t_{2}t_{4}t_{5}t_{6} = t_{11}^{2} \\
t_{2}t_{3}t_{9} = t_{8}t_{12}, \;
t_{2}t_{3}t_{4} = t_{12}^{2}, \;
t_{1}t_{12}t_{13} = t_{8}t_{10} \\
t_{1}t_{11} = t_{7}t_{9}, \;
t_{1}t_{6}t_{12} = t_{8}t_{9}, \;
t_{1}t_{5}t_{12} = t_{7}t_{10} \\
t_{1}t_{4}t_{13} = t_{9}t_{10}, \;
t_{1}t_{4}t_{6} = t_{9}^{2}, \;
t_{1}t_{3}t_{4}t_{5} = t_{10}^{2} \\
t_{1}t_{2}t_{13} = t_{7}t_{8}, \;
t_{1}t_{2}t_{5} = t_{7}^{2}, \;
t_{1}t_{2}t_{3}t_{6} = t_{8}^{2}
 \end{array}
 \right\}
 \subset {\mathbb C}^{13}
 $$

Of course, one may see that the variables $t_{2}$, $t_{3}$, $t_{4}$, $t_{5}$ and $t_{6}$ are uniquely determined with the variable $t_{1}$. Other variables can also be determined in order to get a lower dimensional model.

\section{Appendix}

$P_{3}(t_{1},t_{3},t_{7},t_{11})=4782969 - 4782969 \rho - 19131876 \rho^2 - 4782969 \rho^3 - 4782969 \rho^4 - 4782969 \rho^5 - 1594323 t_{1} - 23914845 \rho t_{1} - 28697814 \rho^2 t_{1} - 1594323 \rho^5 t_{1} - 3188646 t_{1}^2 - 19663317 \rho t_{1}^2 - 6908733 \rho^2 t_{1}^2 + 18068994 \rho^3 t_{1}^2 + 15943230 \rho^4 t_{1}^2 + 3188646 \rho^5 t_{1}^2 - 885735 t_{1}^3 - 885735 \rho t_{1}^3 + 6731586 \rho^2 t_{1}^3 + 20017611 \rho^3 t_{1}^3 + 15943230 \rho^4 t_{1}^3 + 3424842 t_{1}^4 + 5609655 \rho t_{1}^4 + 7144929 \rho^2 t_{1}^4 + 13581270 \rho^3 t_{1}^4 + 7322076 \rho^4 t_{1}^4 + 944784 \rho^5 t_{1}^4 + 3601989 t_{1}^5 + 3109914 \rho t_{1}^5 + 4507407 \rho^2 t_{1}^5 + 5235678 \rho^3 t_{1}^5 + 3149280 \rho^4 t_{1}^5 - 177147 \rho^5 t_{1}^5 + 1003833 t_{1}^6 + 1436859 \rho t_{1}^6 + 997272 \rho^2 t_{1}^6 + 1082565 \rho^3 t_{1}^6 + 747954 \rho^4 t_{1}^6 - 997272 \rho^5 t_{1}^6 + 492075 t_{1}^7 + 658287 \rho t_{1}^7 + 205578 \rho^2 t_{1}^7 + 365229 \rho^3 t_{1}^7 + 96228 \rho^4 t_{1}^7 - 56862 \rho^5 t_{1}^7 + 174231 t_{1}^8 + 143613 \rho t_{1}^8 + 104976 \rho^2 t_{1}^8 + 15309 \rho^3 t_{1}^8 + 45927 \rho^4 t_{1}^8 + 5832 \rho^5 t_{1}^8 + 3159 t_{1}^9 + 24300 \rho t_{1}^9 - 3645 \rho^2 t_{1}^9 - 25272 \rho^3 t_{1}^9 - 10692 \rho^4 t_{1}^9 - 26973 \rho^5 t_{1}^9 + 2106 t_{1}^{10} + 3807 \rho t_{1}^{10} - 1944 \rho^2 t_{1}^{10} - 2349 \rho^3 t_{1}^{10} - 6723 \rho^4 t_{1}^{10} - 1701 \rho^5 t_{1}^{10} + 702 t_{1}^{11} - 675 \rho t_{1}^{11} - 432 \rho^2 t_{1}^{11} - 1107 \rho^3 t_{1}^{11} - 1269 \rho^4 t_{1}^{11} - 1377 \rho^5 t_{1}^{11} + 117 t_{1}^{12} - 99 \rho t_{1}^{12} - 297 \rho^2 t_{1}^{12} - 189 \rho^3 t_{1}^{12} - 117 \rho^4 t_{1}^{12} - 360 \rho^5 t_{1}^{12} + 24 t_{1}^{13} + 15 \rho t_{1}^{13} - 36 \rho^2 t_{1}^{13} - 33 \rho^3 t_{1}^{13} + 15 \rho^4 t_{1}^{13} - 27 \rho^5 t_{1}^{13} - 2 t_{1}^{14} + \rho t_{1}^{14} - 6 \rho^2 t_{1}^{14} - 5 \rho^3 t_{1}^{14} - \rho^4 t_{1}^{14} - \rho^5 t_{1}^{14} - 9565938 \rho^3 t_{11} + 9565938 \rho^4 t_{11} - 3188646 t_{1} t_{11} + 6377292 \rho t_{1} t_{11} + 3188646 \rho^2 t_{1} t_{11} + 9565938 \rho^4 t_{1} t_{11} + 6377292 \rho^5 t_{1} t_{11} - 1062882 t_{1}^2 t_{11} + 3188646 \rho t_{1}^2 t_{11} + 5314410 \rho^2 t_{1}^2 t_{11} + 1062882 \rho^3 t_{1}^2 t_{11} + 3188646 \rho^4 t_{1}^2 t_{11} + 3188646 \rho^5 t_{1}^2 t_{11} - 708588 t_{1}^3 t_{11} + 708588 \rho t_{1}^3 t_{11} + 1771470 \rho^2 t_{1}^3 t_{11} + 354294 \rho^4 t_{1}^3 t_{11} + 354294 \rho^5 t_{1}^3 t_{11} - 236196 t_{1}^4 t_{11} + 236196 \rho t_{1}^4 t_{11} + 590490 \rho^2 t_{1}^4 t_{11} + 118098 \rho^4 t_{1}^4 t_{11} + 118098 \rho^5 t_{1}^4 t_{11} - 78732 t_{1}^5 t_{11} + 78732 \rho t_{1}^5 t_{11} + 196830 \rho^2 t_{1}^5 t_{11} + 39366 \rho^4 t_{1}^5 t_{11} + 39366 \rho^5 t_{1}^5 t_{11} - 26244 t_{1}^6 t_{11} + 26244 \rho t_{1}^6 t_{11} + 65610 \rho^2 t_{1}^6 t_{11} + 13122 \rho^4 t_{1}^6 t_{11} + 13122 \rho^5 t_{1}^6 t_{11} - 8748 t_{1}^7 t_{11} + 8748 \rho t_{1}^7 t_{11} + 21870 \rho^2 t_{1}^7 t_{11} + 4374 \rho^3 t_{1}^7 t_{11} + 4374 \rho^5 t_{1}^7 t_{11} - 1458 t_{1}^8 t_{11} + 5832 \rho^2 t_{1}^8 t_{11} - 2916 \rho^4 t_{1}^8 t_{11} - 1458 \rho^5 t_{1}^8 t_{11} - 486 t_{1}^9 t_{11} - 486 \rho t_{1}^9 t_{11} - 486 \rho^3 t_{1}^9 t_{11} - 972 \rho^4 t_{1}^9 t_{11} - 972 \rho^5 t_{1}^9 t_{11} + 4782969 \rho^3 t_{11}^2 - 4782969 \rho^4 t_{11}^2 + 1594323 t_{1} t_{11}^2 - 3188646 \rho t_{1} t_{11}^2 - 1594323 \rho^2 t_{1} t_{11}^2 - 1594323 \rho^3 t_{1} t_{11}^2 - 3188646 \rho^4 t_{1} t_{11}^2 - 3188646 \rho^5 t_{1} t_{11}^2 - 531441 \rho t_{1}^2 t_{11}^2 - 2125764 \rho^2 t_{1}^2 t_{11}^2 - 531441 \rho^3 t_{1}^2 t_{11}^2 - 531441 \rho^5 t_{1}^2 t_{11}^2 + 177147 t_{1}^3 t_{11}^2 + 177147 \rho t_{1}^3 t_{11}^2 + 177147 \rho^3 t_{1}^3 t_{11}^2 + 354294 \rho^4 t_{1}^3 t_{11}^2 + 354294 \rho^5 t_{1}^3 t_{11}^2 - 19131876 t_{3} - 9565938 \rho t_{3} - 9565938 \rho^2 t_{3} - 9565938 \rho^3 t_{3} - 9565938 \rho^4 t_{3} - 9565938 \rho^5 t_{3} - 12754584 t_{1} t_{3} - 6377292 \rho t_{1} t_{3} + 6377292 \rho^3 t_{1} t_{3} - 9565938 \rho^4 t_{1} t_{3} - 3188646 t_{1}^2 t_{3} - 4251528 \rho t_{1}^2 t_{3} + 5314410 \rho^2 t_{1}^2 t_{3} + 4251528 \rho^3 t_{1}^2 t_{3} - 5314410 \rho^4 t_{1}^2 t_{3} - 4251528 \rho^5 t_{1}^2 t_{3} - 1771470 t_{1}^3 t_{3} - 2125764 \rho t_{1}^3 t_{3} + 1062882 \rho^3 t_{1}^3 t_{3} - 4251528 \rho^4 t_{1}^3 t_{3} - 5314410 \rho^5 t_{1}^3 t_{3} + 236196 t_{1}^4 t_{3} - 708588 \rho t_{1}^4 t_{3} - 354294 \rho^2 t_{1}^4 t_{3} + 708588 \rho^3 t_{1}^4 t_{3} - 1889568 \rho^4 t_{1}^4 t_{3} - 2125764 \rho^5 t_{1}^4 t_{3} + 236196 t_{1}^5 t_{3} - 314928 \rho t_{1}^5 t_{3} - 314928 \rho^2 t_{1}^5 t_{3} + 39366 \rho^3 t_{1}^5 t_{3} - 669222 \rho^4 t_{1}^5 t_{3} - 905418 \rho^5 t_{1}^5 t_{3} + 78732 t_{1}^6 t_{3} - 104976 \rho t_{1}^6 t_{3} - 183708 \rho^2 t_{1}^6 t_{3} - 13122 \rho^3 t_{1}^6 t_{3} - 209952 \rho^4 t_{1}^6 t_{3} - 301806 \rho^5 t_{1}^6 t_{3} + 39366 t_{1}^7 t_{3} - 26244 \rho t_{1}^7 t_{3} - 56862 \rho^2 t_{1}^7 t_{3} + 4374 \rho^3 t_{1}^7 t_{3} - 56862 \rho^4 t_{1}^7 t_{3} - 87480 \rho^5 t_{1}^7 t_{3} + 16038 t_{1}^8 t_{3} - 7290 \rho t_{1}^8 t_{3} - 20412 \rho^2 t_{1}^8 t_{3} - 2916 \rho^3 t_{1}^8 t_{3} - 16038 \rho^4 t_{1}^8 t_{3} - 30618 \rho^5 t_{1}^8 t_{3} + 4860 t_{1}^9 t_{3} - 1458 \rho t_{1}^9 t_{3} - 9234 \rho^2 t_{1}^9 t_{3} - 1944 \rho^3 t_{1}^9 t_{3} - 4374 \rho^4 t_{1}^9 t_{3} - 8262 \rho^5 t_{1}^9 t_{3} + 1944 t_{1}^{10} t_{3} - 162 \rho t_{1}^{10} t_{3} - 2268 \rho^2 t_{1}^{10} t_{3} - 486 \rho^3 t_{1}^{10} t_{3} - 324 \rho^4 t_{1}^{10} t_{3} - 972 \rho^5 t_{1}^{10} t_{3} + 270 t_{1}^{11} t_{3} - 54 \rho t_{1}^{11} t_{3} - 594 \rho^2 t_{1}^{11} t_{3} - 324 \rho^3 t_{1}^{11} t_{3} + 108 \rho^4 t_{1}^{11} t_{3} - 162 \rho^5 t_{1}^{11} t_{3} + 18 t_{1}^{12} t_{3} + 18 \rho t_{1}^{12} t_{3} - 108 \rho^2 t_{1}^{12} t_{3} - 18 \rho^3 t_{1}^{12} t_{3} + 54 \rho^4 t_{1}^{12} t_{3} + 36 \rho^5 t_{1}^{12} t_{3} + 6 t_{1}^{13} t_{3} + 6 \rho t_{1}^{13} t_{3} + 6 \rho^3 t_{1}^{13} t_{3} + 12 \rho^4 t_{1}^{13} t_{3} + 12 \rho^5 t_{1}^{13} t_{3} - 9565938 \rho t_{11} t_{3} + 9565938 \rho^3 t_{11} t_{3} + 3188646 t_{1} t_{11} t_{3} + 3188646 \rho t_{1} t_{11} t_{3} + 6377292 \rho^2 t_{1} t_{11} t_{3} + 6377292 \rho^3 t_{1} t_{11} t_{3} + 6377292 \rho^4 t_{1} t_{11} t_{3} - 3188646 \rho^5 t_{1} t_{11} t_{3} + 6377292 \rho t_{1}^2 t_{11} t_{3} + 5314410 \rho^2 t_{1}^2 t_{11} t_{3} + 1062882 \rho^3 t_{1}^2 t_{11} t_{3} + 1062882 \rho^4 t_{1}^2 t_{11} t_{3} + 1062882 \rho^5 t_{1}^2 t_{11} t_{3} - 1771470 t_{1}^3 t_{11} t_{3} - 708588 \rho t_{1}^3 t_{11} t_{3} + 708588 \rho^2 t_{1}^3 t_{11} t_{3} - 1417176 \rho^3 t_{1}^3 t_{11} t_{3} - 2480058 \rho^4 t_{1}^3 t_{11} t_{3} - 1771470 \rho^5 t_{1}^3 t_{11} t_{3} - 236196 \rho t_{1}^4 t_{11} t_{3} + 118098 \rho^3 t_{1}^4 t_{11} t_{3} - 118098 \rho^4 t_{1}^4 t_{11} t_{3} - 590490 \rho^5 t_{1}^4 t_{11} t_{3} + 39366 t_{1}^5 t_{11} t_{3} + 236196 \rho t_{1}^5 t_{11} t_{3} + 236196 \rho^2 t_{1}^5 t_{11} t_{3} + 78732 \rho^3 t_{1}^5 t_{11} t_{3} + 196830 \rho^4 t_{1}^5 t_{11} t_{3} + 39366 \rho^5 t_{1}^5 t_{11} t_{3} - 26244 t_{1}^6 t_{11} t_{3} + 13122 \rho t_{1}^6 t_{11} t_{3} + 78732 \rho^2 t_{1}^6 t_{11} t_{3} - 13122 \rho^3 t_{1}^6 t_{11} t_{3} - 39366 \rho^4 t_{1}^6 t_{11} t_{3} - 13122 \rho^5 t_{1}^6 t_{11} t_{3} - 8748 t_{1}^7 t_{11} t_{3} - 8748 \rho t_{1}^7 t_{11} t_{3} - 8748 \rho^3 t_{1}^7 t_{11} t_{3} - 17496 \rho^4 t_{1}^7 t_{11} t_{3} - 17496 \rho^5 t_{1}^7 t_{11} t_{3} - 4782969 t_{3}^2 - 4782969 \rho t_{3}^2 - 9565938 \rho^2 t_{3}^2 - 4782969 \rho^3 t_{3}^2 - 4782969 \rho^4 t_{3}^2 - 4782969 \rho^5 t_{3}^2 - 7971615 t_{1} t_{3}^2 - 1594323 \rho t_{1} t_{3}^2 - 4782969 \rho^2 t_{1} t_{3}^2 - 6377292 \rho^3 t_{1} t_{3}^2 - 1594323 \rho^4 t_{1} t_{3}^2 - 8503056 t_{1}^2 t_{3}^2 - 4251528 \rho t_{1}^2 t_{3}^2 - 1594323 \rho^2 t_{1}^2 t_{3}^2 - 2657205 \rho^3 t_{1}^2 t_{3}^2 - 2657205 \rho^4 t_{1}^2 t_{3}^2 + 1062882 \rho^5 t_{1}^2 t_{3}^2 - 2302911 t_{1}^3 t_{3}^2 - 2480058 \rho t_{1}^3 t_{3}^2 - 1417176 \rho^2 t_{1}^3 t_{3}^2 + 177147 \rho^3 t_{1}^3 t_{3}^2 - 354294 \rho^4 t_{1}^3 t_{3}^2 - 1062882 \rho^5 t_{1}^3 t_{3}^2 + 236196 t_{1}^4 t_{3}^2 + 649539 \rho t_{1}^4 t_{3}^2 + 236196 \rho^3 t_{1}^4 t_{3}^2 + 1299078 \rho^4 t_{1}^4 t_{3}^2 + 885735 \rho^5 t_{1}^4 t_{3}^2 - 413343 t_{1}^5 t_{3}^2 - 137781 \rho t_{1}^5 t_{3}^2 - 275562 \rho^3 t_{1}^5 t_{3}^2 + 413343 \rho^5 t_{1}^5 t_{3}^2 - 91854 t_{1}^6 t_{3}^2 - 229635 \rho t_{1}^6 t_{3}^2 - 91854 \rho^2 t_{1}^6 t_{3}^2 - 45927 \rho^4 t_{1}^6 t_{3}^2 + 32805 t_{1}^7 t_{3}^2 + 17496 \rho t_{1}^7 t_{3}^2 - 10935 \rho^2 t_{1}^7 t_{3}^2 + 32805 \rho^3 t_{1}^7 t_{3}^2 + 63423 \rho^4 t_{1}^7 t_{3}^2 + 48114 \rho^5 t_{1}^7 t_{3}^2 - 1458 t_{1}^8 t_{3}^2 + 5832 \rho t_{1}^8 t_{3}^2 + 2187 \rho^2 t_{1}^8 t_{3}^2 - 2187 \rho^3 t_{1}^8 t_{3}^2 + 5832 \rho^4 t_{1}^8 t_{3}^2 + 20412 \rho^5 t_{1}^8 t_{3}^2 - 1215 t_{1}^9 t_{3}^2 - 3159 \rho t_{1}^9 t_{3}^2 - 2673 \rho^2 t_{1}^9 t_{3}^2 - 486 \rho^3 t_{1}^9 t_{3}^2 - 2187 \rho^4 t_{1}^9 t_{3}^2 + 1215 \rho^5 t_{1}^9 t_{3}^2 + 486 t_{1}^{10} t_{3}^2 - 486 \rho^2 t_{1}^{10} t_{3}^2 + 486 \rho^3 t_{1}^{10} t_{3}^2 + 729 \rho^4 t_{1}^{10} t_{3}^2 + 486 \rho^5 t_{1}^{10} t_{3}^2 + 81 t_{1}^{11} t_{3}^2 + 81 \rho t_{1}^{11} t_{3}^2 + 81 \rho^3 t_{1}^{11} t_{3}^2 + 162 \rho^4 t_{1}^{11} t_{3}^2 + 162 \rho^5 t_{1}^{11} t_{3}^2
=0$ 

\s

$
P_{4}(t_{1},t_{3},t_{7},t_{11})=1549681956 + 1937102445 \rho - 387420489 \rho^2 + 1937102445 \rho^4 +  \\
387420489 \rho^5 + 1549681956 t_{1} + 2066242608 \rho t_{1} - 1678822119 \rho^2 t_{1} - 516560652 \rho^3 t_{1} + \\
2711943423 \rho^4 t_{1} + 2195382771 \rho^5 t_{1} + 473513931 t_{1}^2 + 344373768 \rho t_{1}^2 - 645700815 \rho^2 t_{1}^2 + \\
688747536 \rho^3 t_{1}^2 + 4347718821 \rho^4 t_{1}^2 + 2927177028 \rho^5 t_{1}^2 + 301327047 t_{1}^3 + 2424965283 \rho t_{1}^3 \\
+ 1922753538 \rho^2 t_{1}^3 + 3156759540 \rho^3 t_{1}^3 + 5337793404 \rho^4 t_{1}^3 + 2826734679 \rho^5 t_{1}^3 + \\
975725676 t_{1}^4 + 3271550796 \rho t_{1}^4 + 3348078300 \rho^2 t_{1}^4 + 2774122020 \rho^3 t_{1}^4 + 3137627664 \rho^4 t_{1}^4 + \\
2228863554 \rho^5 t_{1}^4 + 435250179 t_{1}^5 + 1551276279 \rho t_{1}^5 + 1729840455 \rho^2 t_{1}^5 + 691936182 \rho^3 t_{1}^5 + \\
680775921 \rho^4 t_{1}^5 + 379448874 \rho^5 t_{1}^5 + 171124002 t_{1}^6 + 628694703 \rho t_{1}^6 + 438970266 \rho^2 t_{1}^6 + \\
48361131 \rho^3 t_{1}^6 + 11160261 \rho^4 t_{1}^6 + 81841914 \rho^5 t_{1}^6 + 17183259 t_{1}^7 + 150752097 \rho t_{1}^7 + \\
173958354 \rho^2 t_{1}^7 - 86802030 \rho^3 t_{1}^7 - 61292862 \rho^4 t_{1}^7 + 105048171 \rho^5 t_{1}^7 - 65071998 t_{1}^8 \\
- 41157153 \rho t_{1}^8 + 16828965 \rho^2 t_{1}^8 - 56568942 \rho^3 t_{1}^8 - 83495286 \rho^4 t_{1}^8 - 30941676 \rho^5 t_{1}^8 \\
- 2775303 t_{1}^9 - 11475189 \rho t_{1}^9 - 3680721 \rho^2 t_{1}^9 + 5708070 \rho^3 t_{1}^9 - 12931731 \rho^4 t_{1}^9 - 7361442 \rho^5 t_{1}^9 + \\ 
2617839 t_{1}^{10} + 3155841 \rho t_{1}^{10} + 5911461 \rho^2 t_{1}^{10} + 3083670 \rho^3 t_{1}^{10} + 5589972 \rho^4 t_{1}^{10} + 5865534 \rho^5 t_{1}^{10} \\
- 2591595 t_{1}^{11} - 756702 \rho t_{1}^{11} + 1791153 \rho^2 t_{1}^{11} + 341172 \rho^3 t_{1}^{11} - 894483 \rho^4 t_{1}^{11} + 150903 \rho^5 t_{1}^{11}\\
 - 127575 t_{1}^{12} - 423549 \rho t_{1}^{12} + 316386 \rho^2 t_{1}^{12} + 898128 \rho^3 t_{1}^{12} - 311283 \rho^4 t_{1}^{12} - 352107 \rho^5 t_{1}^{12} + \\
  219429 t_{1}^{13} + 141183 \rho t_{1}^{13} + 180306 \rho^2 t_{1}^{13} + 275562 \rho^3 t_{1}^{13} + 221130 \rho^4 t_{1}^{13} - 1701 \rho^5 t_{1}^{13} + \\
  15633 t_{1}^{14} + 61641 \rho t_{1}^{14} + 43578 \rho^2 t_{1}^{14} + 31752 \rho^3 t_{1}^{14} + 31590 \rho^4 t_{1}^{14} + 11421 \rho^5 t_{1}^{14} + \\
  1647 t_{1}^{15} + 5157 \rho t_{1}^{15} + 7209 \rho^2 t_{1}^{15} + 4806 \rho^3 t_{1}^{15} - 3402 \rho^4 t_{1}^{15} + 459 \rho^5 t_{1}^{15} + 729 t_{1}^{16} + \\
324 \rho t_{1}^{16} + 495 \rho^2 t_{1}^{16} + 333 \rho^3 t_{1}^{16} - 288 \rho^4 t_{1}^{16} - 774 \rho^5 t_{1}^{16} + 63 t_{1}^{17} + 108 \rho t_{1}^{17} - 18 \rho^2 t_{1}^{17} - 12 \rho^3 t_{1}^{17} \\
+ 3 \rho^4 t_{1}^{17} - 39 \rho^5 t_{1}^{17} + t_{1}^{18} + 5 \rho t_{1}^{18} - \rho^3 t_{1}^{18} - 2 \rho^4 t_{1}^{18} + 4 \rho^5 t_{1}^{18} - 774840978 t_{11} + 774840978 \rho^5 t_{11}\\
 - 1033121304 t_{1} t_{11} - 516560652 \rho t_{1} t_{11} + 1291401630 \rho^2 t_{1} t_{11} + 774840978 \rho^3 t_{1} t_{11} +\\
  516560652 \rho^4 t_{1} t_{11} +  774840978 \rho^5 t_{1} t_{11} - 774840978 t_{1}^2 t_{11} + 258280326 \rho t_{1}^2 t_{11} + \\
  1463588514 \rho^2 t_{1}^2 t_{11} + 1635775398 \rho^3 t_{1}^2 t_{11} + 172186884 \rho^4 t_{1}^2 t_{11} - 344373768 \rho^5 t_{1}^2 t_{11} + \\
  631351908 t_{1}^3 t_{11} + 688747536 \rho t_{1}^3 t_{11} + 1348797258 \rho^2 t_{1}^3 t_{11} + 1234006002 \rho^3 t_{1}^3 t_{11} \\
  - 114791256 \rho^4 t_{1}^3 t_{11} - 373071582 \rho^5 t_{1}^3 t_{11} + 373071582 t_{1}^4 t_{11} + 325241892 \rho t_{1}^4 t_{11} + \\
 296544078 \rho^2 t_{1}^4 t_{11} + 47829690 \rho^3 t_{1}^4 t_{11} - 315675954 \rho^4 t_{1}^4 t_{11} - 660049722 \rho^5 t_{1}^4 t_{11} + \\
 207261990 t_{1}^5 t_{11} + 149866362 \rho t_{1}^5 t_{11} - 117979902 \rho^2 t_{1}^5 t_{11} - 108413964 \rho^3 t_{1}^5 t_{11}\\
  - 255091680 \rho^4 t_{1}^5 t_{11} - 232771158 \rho^5 t_{1}^5 t_{11} + 89282088 t_{1}^6 t_{11} - 1062882 \rho t_{1}^6 t_{11} \\
   - 89282088 \rho^2 t_{1}^6 t_{11} - 127545840 \rho^3 t_{1}^6 t_{11} - 80779032 \rho^4 t_{1}^6 t_{11} - 51018336 \rho^5 t_{1}^6 t_{11} \\
    - 7085880 t_{1}^7 t_{11} - 19840464 \rho t_{1}^7 t_{11} - 62001450 \rho^2 t_{1}^7 t_{11} - 62001450 \rho^3 t_{1}^7 t_{11} \\
  - 24800580 \rho^4 t_{1}^7 t_{11} - 15234642 \rho^5 t_{1}^7 t_{11} - 2834352 t_{1}^8 t_{11} - 7203978 \rho t_{1}^8 t_{11} - 14644152 \rho^2 t_{1}^8 t_{11} \\
  - 10274526 \rho^3 t_{1}^8 t_{11} - 1889568 \rho^4 t_{1}^8 t_{11} + 7085880 \rho^5 t_{1}^8 t_{11} - 1850202 t_{1}^9 t_{11} - 2598156 \rho t_{1}^9 t_{11} \\
  - 1771470 \rho^2 t_{1}^9 t_{11} - 1850202 \rho^3 t_{1}^9 t_{11} + 1299078 \rho^4 t_{1}^9 t_{11} + 1810836 \rho^5 t_{1}^9 t_{11} - 656100 t_{1}^{10} t_{11} \\
  - 314928 \rho t_{1}^{10} t_{11} - 157464 \rho^2 t_{1}^{10} t_{11} + 170586 \rho^3 t_{1}^{10} t_{11} + 328050 \rho^4 t_{1}^{10} t_{11} + 354294 \rho^5 t_{1}^{10} t_{11} \\
  - 17496 t_{1}^{11} t_{11} + 4374 \rho t_{1}^{11} t_{11} + 109350 \rho^2 t_{1}^{11} t_{11} + 131220 \rho^3 t_{1}^{11} t_{11} + 83106 \rho^4 t_{1}^{11} t_{11} + \\
  87480 \rho^5 t_{1}^{11} t_{11}  - 2916 t_{1}^{12} t_{11} + 2916 \rho t_{1}^{12} t_{11} + 13122 \rho^2 t_{1}^{12} t_{11} + 8748 \rho^3 t_{1}^{12} t_{11} + \\
  4374 \rho^4 t_{1}^{12} t_{11} - 5832 \rho^5 t_{1}^{12} t_{11} + 486 \rho t_{1}^{13} t_{11} + 486 \rho^3 t_{1}^{13} t_{11} - 486 \rho^4 t_{1}^{13} t_{11} \\
  - 486 \rho^5 t_{1}^{13} t_{11} - 387420489 \rho t_{11}^2 + 516560652 \rho^3 t_{1} t_{11}^2 + 516560652 \rho^4 t_{1} t_{11}^2 + 430467210 t_{1}^2 t_{11}^2 \\
  + 688747536 \rho t_{1}^2 t_{11}^2 + 688747536 \rho^2 t_{1}^2 t_{11}^2 + 688747536 \rho^3 t_{1}^2 t_{11}^2  + 688747536 \rho^4 t_{1}^2 t_{11}^2 + \\
  430467210 \rho^5 t_{1}^2 t_{11}^2 + 57395628 t_{1}^3 t_{11}^2 + 344373768 \rho t_{1}^3 t_{11}^2 + 344373768 \rho^2 t_{1}^3 t_{11}^2 \\
  + 57395628 \rho^3 t_{1}^3 t_{11}^2 + 4782969 t_{1}^4 t_{11}^2 + 4782969 \rho t_{1}^4 t_{11}^2 - 71744535 \rho^3 t_{1}^4 t_{11}^2 - 167403915 \rho^4 t_{1}^4 t_{11}^2 \\
  - 71744535 \rho^5 t_{1}^4  t_{11}^2 - 31886460 \rho t_{1}^5 t_{11}^2 - 38263752 \rho^2 t_{1}^5 t_{11}^2 - 38263752 \rho^3 t_{1}^5 t_{11}^2 - 38263752 \rho^4 t_{1}^5 t_{11}^2 \\
   - 31886460 \rho^5 t_{1}^5 t_{11}^2 - 3188646 \rho t_{1}^6 t_{11}^2 - 8503056 \rho^2 t_{1}^6 t_{11}^2 - 3188646 \rho^3 t_{1}^6 t_{11}^2 + 708588 \rho^4 t_{1}^7 t_{11}^2 + \\
   708588 \rho^5 t_{1}^7 t_{11}^2 - 59049 t_{1}^8 t_{11}^2 - 1549681956 t_{3} - 3099363912 \rho t_{3} - 2324522934 \rho^2 t_{3} \\
    - 1549681956 \rho^4 t_{3} - 2324522934 \rho^5 t_{3}  + 258280326 t_{1} t_{3} - 1807962282 \rho t_{1} t_{3} \\
    - 774840978 \rho^2 t_{1} t_{3} + 3615924564 \rho^3 t_{1} t_{3} + 2582803260 \rho^4 t_{1} t_{3} - 2066242608 \rho^5 t_{1} t_{3}  +\\
     2496709818 t_{1}^2 t_{3} + 3529831122 \rho t_{1}^2 t_{3} + 3185457354 \rho^2 t_{1}^2 t_{3} + 5423886846 \rho^3 t_{1}^2 t_{3} +\\
    4304672100 \rho^4 t_{1}^2 t_{3}   - 258280326 \rho^5 t_{1}^2 t_{3} +  2955874842 t_{1}^3 t_{3} + 3271550796 \rho t_{1}^3 t_{3} + \\
    3041968284 \rho^2 t_{1}^3 t_{3} + 2152336050 \rho^3 t_{1}^3 t_{3} + 1348797258 \rho^4 t_{1}^3 t_{3} - 1119214746 \rho^5 t_{1}^3 t_{3} +\\
     1788830406 t_{1}^4 t_{3}  + 1769698530 \rho t_{1}^4 t_{3} + 143489070 \rho^2 t_{1}^4 t_{3} - 411335334 \rho^3 t_{1}^4 t_{3} - 28697814 \rho^4 t_{1}^4 t_{3}\\
      - 1320099444 \rho^5 t_{1}^4 t_{3} + 886443588 t_{1}^5 t_{3} + 924707340 \rho t_{1}^5 t_{3} - 197696052 \rho^2 t_{1}^5 t_{3} - 717445350 \rho^3 t_{1}^5 t_{3} \\
      - 31886460 \rho^4 t_{1}^5 t_{3} + 229582512 \rho^5 t_{1}^5 t_{3} - 46766808 t_{1}^6 t_{3} - 100973790 \rho t_{1}^6 t_{3} - 177501294 \rho^2 t_{1}^6 t_{3} \\
       - 640917846 \rho^3 t_{1}^6 t_{3}  - 248714388 \rho^4 t_{1}^6 t_{3} + 61647156 \rho^5 t_{1}^6 t_{3} - 91053558 t_{1}^7 t_{3} - 157306536 \rho t_{1}^7 t_{3} \\
        - 186004350 \rho^2 t_{1}^7 t_{3}  - 161558064 \rho^3 t_{1}^7 t_{3} - 106288200 \rho^4 t_{1}^7 t_{3} - 21966228 \rho^5 t_{1}^7 t_{3} + 8621154 t_{1}^8 t_{3} \\
         - 14526054 \rho t_{1}^8 t_{3}  - 14053662 \rho^2 t_{1}^8 t_{3} + 1535274 \rho^3 t_{1}^8 t_{3} + 46058220 \rho^4 t_{1}^8 t_{3} + 66607272 \rho^5 t_{1}^8 t_{3} \\
         - 18620118 t_{1}^9 t_{3} - 7518906 \rho t_{1}^9 t_{3} + 12282192 \rho^2 t_{1}^9 t_{3} + 3542940 \rho^3 t_{1}^9 t_{3} + 13581270 \rho^4 t_{1}^9 t_{3} + \\
         19053144 \rho^5 t_{1}^9 t_{3}  - 4566456 t_{1}^{10} t_{3} - 2768742 \rho t_{1}^{10} t_{3} + 5038848 \rho^2 t_{1}^{10} t_{3} + 9316620 \rho^3 t_{1}^{10} t_{3} + \\
         1299078 \rho^4 t_{1}^{10} t_{3} + 1141614 \rho^5 t_{1}^{10} t_{3} + 1753974 t_{1}^{11} t_{3} + 1548396 \rho t_{1}^{11} t_{3} + 2934954 \rho^2 t_{1}^{11} t_{3} + \\
         3831624 \rho^3 t_{1}^{11} t_{3} + 2370708 \rho^4 t_{1}^{11} t_{3} + 603612 \rho^5 t_{1}^{11} t_{3} + 237654 t_{1}^{12} t_{3} + 720252 \rho t_{1}^{12} t_{3} + \\
         733374 \rho^2 t_{1}^{12} t_{3} + 542376 \rho^3 t_{1}^{12} t_{3} + 371790 \rho^4 t_{1}^{12} t_{3} + 37908 \rho^5 t_{1}^{12} t_{3} + 45198 t_{1}^{13} t_{3} + \\
          93798 \rho t_{1}^{13} t_{3} + 104976 \rho^2 t_{1}^{13} t_{3} + 78732 \rho^3 t_{1}^{13} t_{3} - 52974 \rho^4 t_{1}^{13} t_{3} - 28188 \rho^5 t_{1}^{13} t_{3} + \\
          18144 t_{1}^{14} t_{3} + 9072 \rho t_{1}^{14} t_{3} + 6156 \rho^2 t_{1}^{14} t_{3} + 2430 \rho^3 t_{1}^{14} t_{3} - 6642 \rho^4 t_{1}^{14} t_{3} - 13284 \rho^5 t_{1}^{14} t_{3} + \\
     1296 t_{1}^{15} t_{3} + 1728 \rho t_{1}^{15} t_{3} - 1350 \rho^2 t_{1}^{15} t_{3} - 1458 \rho^3 t_{1}^{15} t_{3}  - 432 \rho^4 t_{1}^{15} t_{3} - 918 \rho^5 t_{1}^{15} t_{3} + \\
     54 t_{1}^{16} t_{3} + 54 \rho t_{1}^{16} t_{3} - 108 \rho^2 t_{1}^{16} t_{3} - 108 \rho^3 t_{1}^{16} t_{3} - 54 \rho^4 t_{1}^{16} t_{3} + 
    162 \rho^5 t_{1}^{16} t_{3} - 6 \rho t_{1}^{17} t_{3} \\
    - 6 \rho^3 t_{1}^{17} t_{3} + 6 \rho^4 t_{1}^{17} t_{3} + 6 \rho^5 t_{1}^{17} t_{3} + 774840978 \rho^3 t_{11} t_{3} + 774840978 \rho^4 t_{11} t_{3} \\
    - 774840978 \rho^5 t_{11} t_{3} + 1291401630 t_{1} t_{11} t_{3} + 2066242608 \rho t_{1} t_{11} t_{3} + 1807962282 \rho^2 t_{1} t_{11} t_{3} + \\
     2324522934 \rho^3 t_{1} t_{11} t_{3} +1549681956 \rho^4 t_{1} t_{11} t_{3} + 516560652 \rho^5 t_{1} t_{11} t_{3} + 1119214746 t_{1}^2 t_{11} t_{3} + \\
     1980149166 \rho t_{1}^2 t_{11} t_{3} + 2066242608 \rho^2 t_{1}^2 t_{11} t_{3} + 688747536 \rho^3 t_{1}^2 t_{11} t_{3} + 258280326 \rho^4 t_{1}^2 t_{11} t_{3} \\
     - 430467210 \rho^5 t_{1}^2 t_{11} t_{3} + 459165024 t_{1}^3 t_{11} t_{3} + 803538792 \rho t_{1}^3 t_{11} t_{3} + 57395628 \rho^2 t_{1}^3 t_{11} t_{3}\\
      - 602654094 \rho^3 t_{1}^3 t_{11} t_{3} - 1262703816 \rho^4 t_{1}^3 t_{11} t_{3} - 947027862 \rho^5 t_{1}^3 t_{11} t_{3} + 325241892 t_{1}^4 t_{11} t_{3}\\
       - 153055008 \rho t_{1}^4 t_{11} t_{3} - 535692528 \rho^2 t_{1}^4 t_{11} t_{3} - 822670668 \rho^3 t_{1}^4 t_{11} t_{3} - 822670668 \rho^4 t_{1}^4 t_{11} t_{3}\\
        - 420901272 \rho^5 t_{1}^4 t_{11} t_{3} - 70150212 t_{1}^5 t_{11} t_{3} - 267846264 \rho t_{1}^5 t_{11} t_{3} - 491051484 \rho^2 t_{1}^5 t_{11} t_{3} \\
        - 491051484 \rho^3 t_{1}^5 t_{11} t_{3} - 267846264 \rho^4 t_{1}^5 t_{11} t_{3} - 178564176 \rho^5 t_{1}^5 t_{11} t_{3} - 37200870 t_{1}^6 t_{11} t_{3} \\
        - 96722262 \rho t_{1}^6 t_{11} t_{3} - 163683828 \rho^2 t_{1}^6 t_{11} t_{3} - 96722262 \rho^3 t_{1}^6 t_{11} t_{3} - 22320522 \rho^4 t_{1}^6 t_{11} t_{3} + \\
        37200870 \rho^5 t_{1}^6 t_{11} t_{3} - 14880348 t_{1}^7 t_{11} t_{3} - 27280638 \rho t_{1}^7 t_{11} t_{3} - 19840464 \rho^2 t_{1}^7 t_{11} t_{3} \\
        - 10274526 \rho^3 t_{1}^7 t_{11} t_{3} + 24446286 \rho^4 t_{1}^7 t_{11} t_{3} + 25154874 \rho^5 t_{1}^7 t_{11} t_{3} - 7203978 t_{1}^8 t_{11} t_{3} \\
        - 1771470 \rho t_{1}^8 t_{11} t_{3} + 826686 \rho^2 t_{1}^8 t_{11} t_{3} + 4723920 \rho^3 t_{1}^8 t_{11} t_{3} + 8384958 \rho^4 t_{1}^8 t_{11} t_{3} + \\
        6377292 \rho^5 t_{1}^8 t_{11} t_{3} - 236196 t_{1}^9 t_{11} t_{3} + 747954 \rho t_{1}^9 t_{11} t_{3} + 2361960 \rho^2 t_{1}^9 t_{11} t_{3} + 2440692 \rho^3 t_{1}^9 t_{11} t_{3} + \\
        1535274 \rho^4 t_{1}^9 t_{11} t_{3} + 1574640 \rho^5 t_{1}^9 t_{11} t_{3} - 26244 t_{1}^{10} t_{11} t_{3} + 91854 \rho t_{1}^{10} t_{11} t_{3} + 341172 \rho^2 t_{1}^{10} t_{11} t_{3} \\
        + 183708 \rho^3 t_{1}^{10} t_{11} t_{3} + 26244 \rho^4 t_{1}^{10} t_{11} t_{3} - 118098 \rho^5 t_{1}^{10} t_{11} t_{3} + 4374 t_{1}^{11} t_{11} t_{3} + 8748 \rho t_{1}^{11} t_{11} t_{3} + \\
        8748 \rho^3 t_{1}^{11} t_{11} t_{3} - 21870 \rho^4 t_{1}^{11} t_{11} t_{3} - 21870 \rho^5 t_{1}^{11} t_{11} t_{3} + 1458 t_{1}^{12} t_{11} t_{3} +  774840978 t_{11}^2 t_{3} +\\
         774840978 \rho t_{11}^2 t_{3} + 774840978 \rho^2 t_{11}^2 t_{3} + 774840978 \rho^3 t_{11}^2 t_{3} + 774840978 \rho^4 t_{11}^2 t_{3} + \\
          774840978 \rho^5 t_{11}^2 t_{3} + 774840978 \rho t_{1} t_{11}^2 t_{3} + 774840978 \rho^2 t_{1} t_{11}^2 t_{3} - 258280326 \rho^3 t_{1}^2 t_{11}^2 t_{3}\\
           - 774840978 \rho^4 t_{1}^2 t_{11}^2 t_{3} - 258280326 \rho^5 t_{1}^2 t_{11}^2 t_{3} - 229582512 \rho t_{1}^3 t_{11}^2 t_{3} - 258280326 \rho^2 t_{1}^3 t_{11}^2 t_{3} \\
           - 258280326 \rho^3 t_{1}^3 t_{11}^2 t_{3}  - 258280326 \rho^4 t_{1}^3 t_{11}^2 t_{3} - 229582512 \rho^5 t_{1}^3 t_{11}^2 t_{3} - 28697814 \rho t_{1}^4 t_{11}^2 t_{3}\\
            - 86093442 \rho^2 t_{1}^4 t_{11}^2 t_{3} - 28697814 \rho^3 t_{1}^4 t_{11}^2 t_{3} + 9565938 \rho^4 t_{1}^5 t_{11}^2 t_{3} + 9565938 \rho^5 t_{1}^5 t_{11}^2 t_{3} \\
            - 1062882 t_{1}^6 t_{11}^2 t_{3} - 1549681956 t_{3}^2 - 774840978 \rho t_{3}^2 - 387420489 \rho^2 t_{3}^2 + 387420489 \rho^3 t_{3}^2 \\
            - 387420489 \rho^4 t_{3}^2  - 387420489 \rho^5 t_{3}^2 - 1678822119 t_{1} t_{3}^2 - 1420541793 \rho t_{1} t_{3}^2 + 903981141 \rho^2 t_{1} t_{3}^2 +\\
             1033121304 \rho^3 t_{1} t_{3}^2 - 645700815 \rho^4 t_{1} t_{3}^2 - 1549681956 \rho^5 t_{1} t_{3}^2 + 215233605 t_{1}^2 t_{3}^2 - 473513931 \rho t_{1}^2 t_{3}^2 +\\
              903981141 \rho^2 t_{1}^2 t_{3}^2 + 1162261467 \rho^3 t_{1}^2 t_{3}^2 - 516560652 \rho^4 t_{1}^2 t_{3}^2 - 2238429492 \rho^5 t_{1}^2 t_{3}^2 + \\
              1018772397 t_{1}^3 t_{3}^2 + 459165024 \rho t_{1}^3 t_{3}^2 + 57395628 \rho^2 t_{1}^3 t_{3}^2 + 143489070 \rho^3 t_{1}^3 t_{3}^2 \\
          - 258280326 \rho^4 t_{1}^3 t_{3}^2 - 1219657095 \rho^5 t_{1}^3 t_{3}^2 + 511777683 t_{1}^4 t_{3}^2 - 19131876 \rho t_{1}^4 t_{3}^2 - 243931419 \rho^2 t_{1}^4 t_{3}^2 \\
          - 521343621 \rho^3 t_{1}^4 t_{3}^2 - 306110016 \rho^4 t_{1}^4 t_{3}^2 - 325241892 \rho^5 t_{1}^4 t_{3}^2 + 25509168 t_{1}^5 t_{3}^2 - 279006525 \rho t_{1}^5 t_{3}^2 \\
          - 385826166 \rho^2 t_{1}^5 t_{3}^2 - 347562414 \rho^3 t_{1}^5 t_{3}^2 - 172186884 \rho^4 t_{1}^5 t_{3}^2 - 191318760 \rho^5 t_{1}^5 t_{3}^2 + 52612659 t_{1}^6 t_{3}^2 \\
          - 58458510 \rho t_{1}^6 t_{3}^2 - 131265927 \rho^2 t_{1}^6 t_{3}^2 - 23383404 \rho^3 t_{1}^6 t_{3}^2 + 76527504 \rho^4 t_{1}^6 t_{3}^2 + 72807417 \rho^5 t_{1}^6 t_{3}^2 \\
          - 1240029 t_{1}^7 t_{3}^2 - 2657205 \rho t_{1}^7 t_{3}^2 + 14703201 \rho^2 t_{1}^7 t_{3}^2 + 17006112 \rho^3 t_{1}^7 t_{3}^2 + 68378742 \rho^4 t_{1}^7 t_{3}^2 + \\
          65367243 \rho^5 t_{1}^7 t_{3}^2 - 10274526 t_{1}^8 t_{3}^2 + 1299078 \rho t_{1}^8 t_{3}^2 + 16947063 \rho^2 t_{1}^8 t_{3}^2 + 22084326 \rho^3 t_{1}^8 t_{3}^2 + \\
          15234642 \rho^4 t_{1}^8 t_{3}^2 + 11101212 \rho^5 t_{1}^8 t_{3}^2 + 4192479 t_{1}^9 t_{3}^2 + 5727753 \rho t_{1}^9 t_{3}^2 + 10058013 \rho^2 t_{1}^9 t_{3}^2 + \\
          11868849 \rho^3 t_{1}^9 t_{3}^2 + 5747436 \rho^4 t_{1}^9 t_{3}^2 + 2401326 \rho^5 t_{1}^9 t_{3}^2 + 1371249 t_{1}^{10} t_{3}^2 + 2453814 \rho t_{1}^{10} t_{3}^2 + \\
          2775303 \rho^2 t_{1}^{10} t_{3}^2 + 1863324 \rho^3 t_{1}^{10} t_{3}^2 + 807003 \rho^4 t_{1}^{10} t_{3}^2 - 360855 \rho^5 t_{1}^{10} t_{3}^2 + 301806 t_{1}^{11} t_{3}^2 + \\
          422091 \rho t_{1}^{11} t_{3}^2 + 279936 \rho^2 t_{1}^{11} t_{3}^2 + 146529 \rho^3 t_{1}^{11} t_{3}^2 - 334611 \rho^4 t_{1}^{11} t_{3}^2 - 295245 \rho^5 t_{1}^{11} t_{3}^2 + \\
          105705 t_{1}^{12} t_{3}^2 + 35721 \rho t_{1}^{12} t_{3}^2 - 22599 \rho^2 t_{1}^{12} t_{3}^2 - 45198 \rho^3 t_{1}^{12} t_{3}^2 - 64152 \rho^4 t_{1}^{12} t_{3}^2 - 70713 \rho^5 t_{1}^{12} t_{3}^2 +\\
           6561 t_{1}^{13} t_{3}^2 + 2916 \rho t_{1}^{13} t_{3}^2 - 19926 \rho^2 t_{1}^{13} t_{3}^2 - 19926 \rho^3 t_{1}^{13} t_{3}^2 - 6075 \rho^4 t_{1}^{13} t_{3}^2 - 6075 \rho^5 t_{1}^{13} t_{3}^2 + \\
           324 t_{1}^{14} t_{3}^2 - 324 \rho t_{1}^{14} t_{3}^2 - 2106 \rho^2 t_{1}^{14} t_{3}^2 - 1620 \rho^3 t_{1}^{14} t_{3}^2 + 1944 \rho^5 t_{1}^{14} t_{3}^2 - 54 t_{1}^{15} t_{3}^2 \\
           - 108 \rho t_{1}^{15} t_{3}^2 - 108 \rho^3 t_{1}^{15} t_{3}^2 + 162 \rho^4 t_{1}^{15} t_{3}^2 + 162 \rho^5 t_{1}^{15} t_{3}^2 - 9 t_{1}^{16} t_{3}^2 + 774840978 \rho^2 t_{11} t_{3}^2 + \\
           774840978 \rho^4 t_{11} t_{3}^2 - 258280326 t_{1} t_{11} t_{3}^2 + 516560652 \rho t_{1} t_{11} t_{3}^2 + 516560652 \rho^2 t_{1} t_{11} t_{3}^2 +\\
            258280326 \rho^3 t_{1} t_{11} t_{3}^2 - 258280326 \rho^4 t_{1} t_{11} t_{3}^2 - 516560652 \rho^5 t_{1} t_{11} t_{3}^2 + 516560652 t_{1}^2 t_{11} t_{3}^2 \\
            - 258280326 \rho^3 t_{1}^2 t_{11} t_{3}^2 - 602654094 \rho^4 t_{1}^2 t_{11} t_{3}^2 - 344373768 \rho^5 t_{1}^2 t_{11} t_{3}^2 - 315675954 \rho t_{1}^3 t_{11} t_{3}^2 \\
             - 516560652 \rho^2 t_{1}^3 t_{11} t_{3}^2 - 545258466 \rho^3 t_{1}^3 t_{11} t_{3}^2 - 200884698 \rho^4 t_{1}^3 t_{11} t_{3}^2 - 286978140 \rho^5 t_{1}^3 t_{11} t_{3}^2 \\
             - 38263752 t_{1}^4 t_{11} t_{3}^2  - 86093442 \rho t_{1}^4 t_{11} t_{3}^2 - 210450636 \rho^2 t_{1}^4 t_{11} t_{3}^2 - 76527504 \rho^3 t_{1}^4 t_{11} t_{3}^2 + \\
             47829690 \rho^4 t_{1}^4 t_{11} t_{3}^2 + 95659380 \rho^5 t_{1}^4 t_{11} t_{3}^2 - 9565938 t_{1}^5 t_{11} t_{3}^2 - 22320522 \rho t_{1}^5 t_{11} t_{3}^2 + \\
             6377292 \rho^2 t_{1}^5 t_{11} t_{3}^2 + 9565938 \rho^3 t_{1}^5 t_{11} t_{3}^2 + 76527504 \rho^4 t_{1}^5 t_{11} t_{3}^2 + 73338858 \rho^5 t_{1}^5 t_{11} t_{3}^2 \\
             - 17006112 t_{1}^6 t_{11} t_{3}^2 + 2125764 \rho t_{1}^6 t_{11} t_{3}^2 + 12754584 \rho^2 t_{1}^6 t_{11} t_{3}^2 + 14880348 \rho^3 t_{1}^6 t_{11} t_{3}^2 +\\
              23383404 \rho^4 t_{1}^6 t_{11} t_{3}^2 + 15943230 \rho^5 t_{1}^6 t_{11} t_{3}^2 - 354294 t_{1}^7 t_{11} t_{3}^2 + 3897234 \rho t_{1}^7 t_{11} t_{3}^2 +\\
           9920232 \rho^2 t_{1}^7 t_{11} t_{3}^2 + 8857350 \rho^3 t_{1}^7 t_{11} t_{3}^2 + 4251528 \rho^4 t_{1}^7 t_{11} t_{3}^2 + 4960116 \rho^5 t_{1}^7 t_{11} t_{3}^2 + \\
           118098 t_{1}^8 t_{11} t_{3}^2 +472392 \rho t_{1}^8 t_{11} t_{3}^2 + 1889568 \rho^2 t_{1}^8 t_{11} t_{3}^2 + 826686 \rho^3 t_{1}^8 t_{11} t_{3}^2 \\
           - 236196 \rho^4 t_{1}^8 t_{11} t_{3}^2 - 708588 \rho^5 t_{1}^8 t_{11} t_{3}^2 +78732 t_{1}^9 t_{11} t_{3}^2 + 39366 \rho t_{1}^9 t_{11} t_{3}^2 + 39366 \rho^3 t_{1}^9 t_{11} t_{3}^2 \\
            - 196830 \rho^4 t_{1}^9 t_{11} t_{3}^2 - 196830 \rho^5 t_{1}^9 t_{11} t_{3}^2 + 26244 t_{1}^{10} t_{11} t_{3}^2 - 387420489 \rho^4 t_{11}^2 t_{3}^2 \\
             - 258280326 \rho t_{1} t_{11}^2 t_{3}^2 - 258280326 \rho^2 t_{1} t_{11}^2 t_{3}^2 - 258280326 \rho^3 t_{1} t_{11}^2 t_{3}^2 - 258280326 \rho^4 t_{1} t_{11}^2 t_{3}^2 \\
              - 258280326 \rho^5 t_{1} t_{11}^2 t_{3}^2 - 43046721 \rho t_{1}^2 t_{11}^2 t_{3}^2 - 172186884 \rho^2 t_{1}^2 t_{11}^2 t_{3}^2 - 43046721 \rho^3 t_{1}^2 t_{11}^2 t_{3}^2 \\
              + 28697814 \rho^4 t_{1}^3 t_{11}^2 t_{3}^2 + 28697814 \rho^5 t_{1}^3 t_{11}^2 t_{3}^2 - 4782969 t_{1}^4 t_{11}^2 t_{3}^2 - 774840978 t_{3}^3 \\
             - 258280326 t_{1} t_{3}^3 + 258280326 \rho t_{1} t_{3}^3 + 774840978 \rho^2 t_{1} t_{3}^3 + 774840978 \rho^3 t_{1} t_{3}^3 + 258280326 \rho^5 t_{1} t_{3}^3 \\
             + 172186884 \rho t_{1}^2 t_{3}^3 + 602654094 \rho^2 t_{1}^2 t_{3}^3 + 344373768 \rho^3 t_{1}^2 t_{3}^3 - 172186884 \rho^4 t_{1}^2 t_{3}^3 \\
              - 344373768 \rho^5 t_{1}^2 t_{3}^3 + 57395628 t_{1}^3 t_{3}^3 + 114791256 \rho t_{1}^3 t_{3}^3 + 114791256 \rho^2 t_{1}^3 t_{3}^3 + \\
              86093442 \rho^3 t_{1}^3 t_{3}^3 - 258280326 \rho^4 t_{1}^3 t_{3}^3 - 315675954 \rho^5 t_{1}^3 t_{3}^3 + 86093442 t_{1}^4 t_{3}^3 + \\
              38263752 \rho t_{1}^4 t_{3}^3 + 9565938 \rho^2 t_{1}^4 t_{3}^3 - 9565938 \rho^3 t_{1}^4 t_{3}^3  - 124357194 \rho^4 t_{1}^4 t_{3}^3 - 133923132 \rho^5 t_{1}^4 t_{3}^3 +\\
               28697814 t_{1}^5 t_{3}^3 + 6377292 \rho t_{1}^5 t_{3}^3 - 25509168 \rho^2 t_{1}^5 t_{3}^3  - 31886460 \rho^3 t_{1}^5 t_{3}^3 - 51018336 \rho^4 t_{1}^5 t_{3}^3 \\
               - 60584274 \rho^5 t_{1}^5 t_{3}^3 + 9565938 t_{1}^6 t_{3}^3 + 1062882 \rho t_{1}^6 t_{3}^3 - 14880348 \rho^2 t_{1}^6 t_{3}^3 - 13817466 \rho^3 t_{1}^6 t_{3}^3 \\
                - 17006112 \rho^4 t_{1}^6 t_{3}^3 - 17006112 \rho^5 t_{1}^6 t_{3}^3 + 3188646 t_{1}^7 t_{3}^3 - 4960116 \rho^2 t_{1}^7 t_{3}^3 - 4960116 \rho^3 t_{1}^7 t_{3}^3 \\
                - 4960116 \rho^4 t_{1}^7 t_{3}^3 - 4960116 \rho^5 t_{1}^7 t_{3}^3 + 944784 t_{1}^8 t_{3}^3 - 118098 \rho t_{1}^8 t_{3}^3 - 2007666 \rho^2 t_{1}^8 t_{3}^3 \\
                - 2007666 \rho^3 t_{1}^8 t_{3}^3 - 1653372 \rho^4 t_{1}^8 t_{3}^3 - 1771470 \rho^5 t_{1}^8 t_{3}^3 +  275562 t_{1}^9 t_{3}^3 - 78732 \rho t_{1}^9 t_{3}^3 \\
                 - 826686 \rho^2 t_{1}^9 t_{3}^3 - 708588 \rho^3 t_{1}^9 t_{3}^3 - 472392 \rho^4 t_{1}^9 t_{3}^3 - 393660 \rho^5 t_{1}^9 t_{3}^3 + 65610 t_{1}^{10} t_{3}^3 \\
                 - 52488 \rho t_{1}^{10} t_{3}^3 - 236196 \rho^2 t_{1}^{10} t_{3}^3 - 223074 \rho^3 t_{1}^{10} t_{3}^3 - 65610 \rho^4 t_{1}^{10} t_{3}^3 - 39366 \rho^5 t_{1}^{10} t_{3}^3 \\
             - 8748 t_{1}^{11} t_{3}^3 - 17496 \rho t_{1}^{11} t_{3}^3 - 65610 \rho^2 t_{1}^{11} t_{3}^3 - 56862 \rho^3 t_{1}^{11} t_{3}^3 - 4374 \rho^4 t_{1}^{11} t_{3}^3 - 2916 t_{1}^{12} t_{3}^3 \\
              - 2916 \rho t_{1}^{12} t_{3}^3  - 8748 \rho^2 t_{1}^{12} t_{3}^3 - 5832 \rho^3 t_{1}^{12} t_{3}^3 + 2916 \rho^4 t_{1}^{12} t_{3}^3 + 7290 \rho^5 t_{1}^{12} t_{3}^3 \\
              - 972 t_{1}^{13} t_{3}^3 - 486 \rho t_{1}^{13} t_{3}^3 - 486 \rho^3 t_{1}^{13} t_{3}^3 + 972 \rho^4 t_{1}^{13} t_{3}^3 + 972 \rho^5 t_{1}^{13} t_{3}^3 - 162 t_{1}^{14} t_{3}^3 \\
               + 774840978 \rho^2 t_{11} t_{3}^3 + 258280326 \rho^2 t_{1} t_{11} t_{3}^3 - 258280326 \rho^4 t_{1} t_{11} t_{3}^3 - 258280326 \rho^5 t_{1} t_{11} t_{3}^3 + \\
               86093442 t_{1}^2 t_{11} t_{3}^3 + 86093442 \rho^2 t_{1}^2 t_{11} t_{3}^3 - 86093442 \rho^4 t_{1}^2 t_{11} t_{3}^3  - 86093442 \rho^5 t_{1}^2 t_{11} t_{3}^3 + \\
                28697814 t_{1}^3 t_{11} t_{3}^3 + 28697814 \rho^2 t_{1}^3 t_{11} t_{3}^3  - 28697814 \rho^4 t_{1}^3 t_{11} t_{3}^3 - 28697814 \rho^5 t_{1}^3 t_{11} t_{3}^3 + \\
                9565938 t_{1}^4 t_{11} t_{3}^3 + 9565938 \rho^2 t_{1}^4 t_{11} t_{3}^3 - 9565938 \rho^4 t_{1}^4 t_{11} t_{3}^3 - 9565938 \rho^5 t_{1}^4 t_{11} t_{3}^3 + \\
                3188646 t_{1}^5 t_{11} t_{3}^3 + 3188646 \rho^2 t_{1}^5 t_{11} t_{3}^3  - 3188646 \rho^4 t_{1}^5 t_{11} t_{3}^3 - 3188646 \rho^5 t_{1}^5 t_{11} t_{3}^3 + \\
                1062882 t_{1}^6 t_{11} t_{3}^3 + 1062882 \rho^2 t_{1}^6 t_{11} t_{3}^3 - 1062882 \rho^4 t_{1}^6 t_{11} t_{3}^3 - 1062882 \rho^5 t_{1}^6 t_{11} t_{3}^3 + \\
                354294 t_{1}^7 t_{11} t_{3}^3 - 354294 \rho^4 t_{1}^7 t_{11} t_{3}^3 - 354294 \rho^5 t_{1}^7 t_{11} t_{3}^3 +118098 t_{1}^8 t_{11} t_{3}^3 - 387420489 t_{3}^4 \\
                 - 258280326 t_{1} t_{3}^4 - 129140163 t_{1}^2 t_{3}^4 - 57395628 t_{1}^3 t_{3}^4 - 23914845 t_{1}^4 t_{3}^4 - 9565938 t_{1}^5 t_{3}^4 \\            
              - 3720087 t_{1}^6 t_{3}^4 - 1062882 t_{1}^7 t_{3}^4 - 295245 t_{1}^8 t_{3}^4 - 78732 t_{1}^9 t_{3}^4 - 19683 t_{1}^{10} t_{3}^4 \\ 
                 - 4374 t_{1}^{11} t_{3}^4 - 729 t_{1}^{12} t_{3}^4
=0.
$


\end{document}